\documentclass[12pt]{amsart}
\usepackage{amssymb,latexsym}

\numberwithin{equation}{section}

\newtheorem{theorem}{Theorem}[section]
\newtheorem{proposition}[theorem]{Proposition}

\newtheorem{corollary}[theorem]{Corollary}
\newtheorem{lemma}[theorem]{Lemma}

\theoremstyle{definition}

\newtheorem{remark}[theorem]{Remark}
\newtheorem{example}[theorem]{Example}

\def\AA{\mathcal{A}}
\def\BB{\mathbf{B}}
\def\FF{\mathcal{F}}
\def\ZZ{\mathbb{Z}}

\def\mm{\mathbf{m}}

\def\TT{\mathcal{T}}

\def\ZZ{\mathbb{Z}}
\def\PP{\mathbb{P}}
\def\XX{\mathbf{X}}

\begin{document}

\title{Triangular bases in quantum cluster algebras}

\author{Arkady Berenstein}
\address{\noindent Department of Mathematics, University of Oregon,
Eugene, OR 97403, USA} \email{arkadiy@math.uoregon.edu}

\author{Andrei Zelevinsky}
\address{\noindent Department of Mathematics, Northeastern University,
  Boston, MA 02115, USA}
\email{andrei@neu.edu}

\date{June 15, 2012; revised November 11, 2012}

\thanks{Research supported in part by NSF grants DMS-0800247, DMS-1101507 (A.~B.)
and DMS-0801187, DMS-1103813 (A.~Z.)}

\maketitle

\section{Introduction and main results}
\label{sec:intro}

One of the main motivations behind the theory of cluster algebras was to create an algebraic formalism for a better understanding of canonical bases in quantum groups.
To this end a lot of recent activity has been directed towards various constructions of ``natural" bases in cluster algebras.
An overview of these approaches with relevant references can be found in \cite{plamondon}.

In this paper we develop a new approach.
It is in fact much closer to Lusztig's original way of constructing a canonical basis \cite{lusztig}
(and the pioneering construction of the Kazhdan-Lusztig basis in a Hecke algebra).
The key ingredient of our approach is a version of \emph{Lusztig's Lemma} generalizing \cite[Theorem~1.2]{du} (see also \cite{lusztig-Hecke}).
Here is this version (a proof will be given in Section~\ref{sec:proof-Lusztig-lemma}).

\begin{theorem}
\label{th:Lusztig-lemma}
Let $\AA$ be a free $\ZZ[v, v^{-1}]$-module with a basis $\{E_a: a\in
L\}$ indexed by a partially ordered set $(L, \prec)$ such that,
for any $a \in L$, the lengths of chains in~$L$ with the top element~$a$ are bounded from above.
Let $x \mapsto \overline x$ be a $\ZZ$-linear involution on $\AA$ such that,
for all $f \in \ZZ[v, v^{-1}]$ and $x\in \AA$, we have
\begin{equation}
\label{eq:bar-v-to-v-inverse}
\text{$\overline {f\ x}=\overline f \ \overline x$, where $\overline f(v) = f(v^{-1})$.}
\end{equation}
Suppose that
\begin{equation}
\label{eq:bar-triangularity}
\overline {E_a} - E_a \in \oplus_{a' \prec a} \ZZ[v, v^{-1}]
E_{a'} \quad (a \in L).
\end{equation}
Then for every $a\in L$ there exists a unique element $C_a\in \AA$ such that:
\begin{equation}
\label{eq:b-bar-fixed}
\overline {C_a} = C_a;
\end{equation}
\begin{equation}
\label{eq:b-triangularity}
C_a - E_a \in \oplus_{a' \in L} v \ZZ[v]
E_{a'}.
\end{equation}
Moreover, the element $C_a$ satisfies
\begin{equation}
\label{eq:b-double-triangularity}
C_a - E_a \in \oplus_{a' \prec a} v \ZZ[v] E_{a'},
\end{equation}
hence the elements $C_a$ for $a \in L$ form a $\ZZ[v, v^{-1}]$-basis in~$\AA$.
\end{theorem}

 We will apply Theorem~\ref{th:Lusztig-lemma} in the situation where $\AA$ is a \emph{quantum cluster algebra}
(in the sense of \cite{bz-advances}) with an acyclic quantum seed.
To state the main results we need to recall some terminology and notation from
\cite{bz-advances}.

Recall that a (labeled) \emph{quantum seed} is specified by the following data:
\begin{itemize}
\item Two positive integers $m \geq n$.
\item
An $m \! \times \! n$ integer matrix $\tilde B$.
We represent $\tilde B$ as the family of its columns $b_j \in \ZZ^m$ for $j \in [1,n] = \{1, \dots, n\}$.
\item
A family of positive integers $d_1, \dots, d_n$.
\item
A skew-symmetric bilinear form $\Lambda: \ZZ^m\times \ZZ^m\to \ZZ$  satisfying the \emph{compatibility condition} with $\tilde B$:
\begin{equation}
\label{eq:orthogonality}
\Lambda(b_j, e_i)=\delta_{ij} d_j
\end{equation}
for $i\in [1,m]$, $j\in [1,n]$, where $e_1, \dots, e_m$ are the standard basis vectors in $\ZZ^m$.
We identify $\Lambda$ with the skew-symmetric $m \times m$ matrix $\Lambda = (\lambda_{ij})$, where
$\lambda_{ij} = \Lambda(e_i,e_j)$.
\item
The \emph{based quantum torus}
$\TT = \TT(\Lambda)$, that is, the $\ZZ[v, v^{-1}]$-algebra
with a distinguished $\ZZ[v, v^{-1}]$-basis $\{X^e: e \in \ZZ^m\}$ and the
multiplication given by
\begin{equation}
\label{eq:multiplication-quantum-torus}
X^e X^f = v^{\Lambda(e,f)} X^{e+f} \quad (e,f \in \ZZ^m) \, .
\end{equation}
\end{itemize}

We abbreviate $X^{e_i} = X_i$ for $i \in [1,m]$, and call the set $\tilde \XX = \{X_1, \dots, X_m\}$
the \emph{extended cluster} of a quantum seed.
With some abuse of notation, we will denote the quantum seed simply as $(\tilde \XX, \tilde B)$.
We call the $n \times n$ submatrix $B = (b_{ij})_{i, j \in [1,n]}$ of $\tilde B$ the \emph{exchange matrix} of a quantum seed.
Note that the condition that the form $\Lambda$ is skew-symmetric, together with \eqref{eq:orthogonality}, implies that, for every
$j, k \in [1,n]$, we have
\begin{equation}
\label{eq:skew-symm}
d_j b_{jk} = \Lambda(b_j, b_k) = - \Lambda(b_k, b_j) = - d_k b_{kj},
\end{equation}
that is the matrix $B$ is \emph{skew-symmetrizable}.

The elements of $\tilde \XX$ and their inverses generate $\TT$ as a $\ZZ[v, v^{-1}]$-algebra, subject to the
\emph{quasi-commutation relations}
\begin{equation}
\label{eq:Xi-q-com}
X_i X_j = v^{2 \lambda_{ij}} X_j X_i \quad (i, j \in [1,m]) \, .
\end{equation}
We call the set $\XX = \{X_j: j \in [1,n]\}$ the
\emph{cluster} of a quantum seed, and the elements $X_i \in \tilde \XX - \XX$ the \emph{frozen variables}.
Let $\PP$ denote the multiplicative subgroup in $\TT$ generated by $v$ and all frozen variables, and let $\ZZ \PP$ be the integer group ring
of $\PP$, that is, the ring of (non-commutative) Laurent polynomials in the frozen variables with coefficients in $\ZZ[v, v^{-1}]$.
Note that $\TT$ satisfies the \emph{Ore condition} (see e.g. Appendix to \cite{bz-advances}), so can be viewed as a $\ZZ \PP$-subalgebra of the ambient
\emph{skew-field of fractions} $\FF$.

The following example shows that every skew-symmetrizable $n \times n$ matrix $B$
is an exchange matrix of a special quantum seed that we call \emph{principal}.
This construction will play an important role later in the paper.

\begin{example}
{\bf Principal quantization.}
\label{ex:principal-q-seed}
Let $B$ be a skew-symmetrizable integer $n \times n$
matrix, that is, $B$ satisfies \eqref{eq:skew-symm} for some positive
integers $d_1, \dots, d_n$; in other words,  $DB$ is a skew-symmetric matrix, where
$D$ is the diagonal matrix with diagonal entries  $d_1, \dots, d_n$.
We set $m= 2n$, and define a $m \times n$ integer matrix $\tilde B$ as
\begin{equation}
\label{tilde-B-principal}
\tilde B =
\left(\!\!\begin{array}{c}
B \\ I_n \\
\end{array}\!\!\right) \, ,
\end{equation}
where $I_n$ is the identity $n \times n$ matrix.

A direct inspection shows that the bilinear form~$\Lambda$ on $\ZZ^{2n}$
with the matrix
\begin{equation}
\label{eq:lambda-principal}
\Lambda =
\left(\!\!\begin{array}{cc}
0 & -D \\
D & -DB \\
\end{array}\!\!\right)
\end{equation}
is skew-symmetric, and  satisfies \eqref{eq:orthogonality}.
The corresponding extended cluster $\tilde \XX$ consists of cluster variables $X_1, \dots, X_n$
and frozen variables $X_{n+1}, \dots, X_{2n}$.
The quasi-commutation relations \eqref{eq:Xi-q-com} are as follows: all cluster variables
$X_1, \dots, X_n$ commute with each other, and we have
\begin{equation}
\label{eq:Xi-q-com-principal}
X_i X_k = v^{2 d_i \delta_{i,k+n}} X_k X_i, \,\,
X_i X_j = v^{2 d_j b_{ji}} X_j X_i \,\, (k \in [1,n], \,\, i, j \in [n+1,2n]) \, .
\end{equation}
\end{example}

Returning to the general case, for $k \in [1,n]$, the \emph{quantum seed mutation} $\mu_k$
transforms $(\tilde \XX, \tilde B)$ into a quantum seed $(\tilde \XX', \tilde B')$ defined as follows:
\begin{itemize}
\item
The matrix entries $b'_{ij}$ of $\tilde B'$ are given by
\begin{equation}
\label{eq:matrix-mutation}
b'_{ij} =
\left\{\!
\begin{array}{ll}
-b_{ij} & {\rm if~} i=k {\rm~or~} j=k; \\[.05in]
b_{ij} + [b_{ik}]_+ [b_{kj}]_+  - [-b_{ik}]_+ [-b_{kj}]_+ & {\rm otherwise}.
\end{array}
\right.
\end{equation}
\item The extended cluster $\tilde {\XX'}$ is obtained from $\tilde \XX = \{X_1, \dots, X_m\}$ by replacing $X_k$ with $X'_k$ given by
\begin{equation}
\label{eq:q-exchange-rel-xx}
X'_k = X^{- e_k + [b_k]_+} \ + \
 X^{- e_k + [-b_k]_+} \ .
\end{equation}
\end{itemize}
Here, as before, $b_k \in \ZZ^m$ is the $k$th column of $\tilde B$, and, for an integer vector $b$,
we use the notation $[b]_+$ for the vector obtained by applying the operation
$c \mapsto \max(c,0)$ to each component of~$b$ (in particular, $[b]_+ = \max(b,0)$ for $b \in \ZZ$).

Note that $\tilde {\XX'}$ generates a based quantum torus $\TT'$ defined in the same way as the quantum torus $\TT$ above, but with the quasi-commutation relations
governed by the new skew-symmetric form $\Lambda'$ on $\ZZ^m$ given as follows: $\Lambda'(e_i, e_j) = \Lambda(e_i, e_j)$ for $i, j \neq k$, and
$\Lambda'(e_k, e_j) = \Lambda (e'_k, e_j)$, where we abbreviate
\begin{equation}
\label{eq:e-prime-k}
e'_k = - e_k + [b_k]_+ \in \ZZ^m \, .
\end{equation}
Note also that the skew-field of fractions of $\TT' = \TT(\Lambda')$ is naturally identified with~$\FF$.

Now the \emph{quantum cluster algebra} $\AA = \AA(\tilde \XX, \tilde B)$ is defined as the $\ZZ \PP$-subalgebra of the ambient skew-field~$\FF$
generated by all \emph{cluster varables}, that is, by the union of clusters of all quantum seeds obtained from the initial quantum seed by a finite sequence of mutations.
According to the \emph{quantum Laurent phenomenon} proved in \cite{bz-advances}, $\AA$ is a subalgebra in the quantum based torus associated with any of its seeds.

As shown in \cite[Section~7]{bz-advances}, the structure of $\AA(\tilde \XX, \tilde B)$ simplifies dramatically if the initial quantum seed is \emph{acyclic}.
(Recall that this means the following: the directed graph $\Gamma(B)$ with the set of vertices $[1,n]$ and oriented edges
$j \to i$ for all pairs $(i, j)$ such that $b_{ij} > 0$, has no oriented cycles.)
Namely, we have:
\begin{align}
\label{eq:lower-qca}
&\text{For an acyclic quantum seed $(\tilde \XX, \tilde B)$, the cluster algebra $\AA$ is}\\
\nonumber
&\text{a $\ZZ \PP$-subalgebra of $\TT$ generated by the elements $X_k, X'_k$ for $k \in [1,n]$.}
\end{align}

In fact, the results in \cite[Section~7]{bz-advances} provide the following sharper statement.
Let $i \lhd j$ be any linear order on $[1,n]$.
For any $a \in \ZZ^m$, we denote by $E_a^\circ$ the element of $\AA$ given by
\begin{equation}
\label{eq:standard-monomial-non-normalized}
E_a^\circ = X^{\sum_{i \in [1,m] - [1,n]} a_i e_i + \sum_{j \in [1,n]} [a_j]_+ e_j} \cdot
\prod_{k \in [1,n]}^\lhd (X'_k)^{[- a_k]_+} \ ,
\end{equation}
where the notation $\prod_{k \in [1,n]}^\lhd$ means that the product is taken in increasing order with respect to $\lhd$.
Then we have (see \cite[Theorem~7.3 and (7.4)]{bz-advances}):
\begin{align}
\label{eq:monomial-basis-lower-qca}
&\text{For an acyclic quantum seed $(\tilde \XX, \tilde B)$, the elements $E_a^\circ$}\\
\nonumber
&\text{for $a \in \ZZ^m$ form a basis in $\AA$ as a $\ZZ[v, v^{-1}]$-module.}
\end{align}

Let $X \to \overline X$ be the involution of $\TT$ (``bar-involution") satisfying \eqref{eq:bar-v-to-v-inverse}, and such that
$\overline{X^a} = X^a$ for all $a \in \ZZ^m$.
An easy check shows that
\begin{equation}
\label{eq:bar-properties}
\overline{X Y} = \overline{Y} \ \overline{X} \quad (X, Y \in \TT), \, \overline{X'_k} = X'_k \quad (k \in [1,n]) \ .
\end{equation}
It follows that the subalgebra $\AA \subset \TT$ is invariant under the bar-involution.
Furthermore, as shown in \cite{bz-advances}, the bar-involution on $\AA$ thus obtained is independent of the choice of
an initial quantum seed.

For $a \in \ZZ^m$, we define
\begin{equation}
\label{eq:r-of-a}
r(a) = \sum_{k \in [1,n]} [-a_k]_+ \ .
\end{equation}
We make $\ZZ^m$ into a partially ordered set by setting
\begin{equation}
\label{eq:partial-order-on-Zm}
a' \prec a \Longleftrightarrow r(a') < r(a) \, ;
\end{equation}
note that $L = \ZZ^m$ satisfies the boundedness condition in Theorem~\ref{th:Lusztig-lemma}.

We define the \emph{leading term} $LT(E_a^\circ)$ as an element of $\TT$ obtained from the expression \eqref{eq:standard-monomial-non-normalized}
for  $E_a^\circ$ by replacing each $X'_k$ with $X^{e'_k}$.
Then we set
\begin{equation}
\label{eq:E-a}
E_a = v^{\nu(a)} E_a^\circ \, ,
\end{equation}
where the exponent $\nu(a) \in \ZZ$ is determined from the condition that
\begin{equation}
\label{eq:nu-of-a}
\text{the element $v^{\nu(a)} LT(E_a^\circ) \in \TT$ is invariant under the bar-involution.}
\end{equation}

\begin{remark}
\label{rem:leading-term-variant}
It is easy to give an explicit formula for $\nu(a)$ but we will not need it.
Note also that $\nu(a)$ and hence the element $E_a$ will not change if we modify the definition of $LT(E_a^\circ)$
by replacing each $X'_k$ with $X^{e'_k - b_k} = X^{- e_k + [-b_k]_+}$ instead of $X^{e'_k}$.
\end{remark}

Finally everything is ready for stating our application of Theorem~\ref{th:Lusztig-lemma}.

\begin{theorem}
\label{th:Lusztig-lemma-ca}
Let $\AA = \AA(\tilde \XX, \tilde B)$ be the cluster algebra associated with an acyclic quantum seed $(\tilde \XX, \tilde B)$, and
$\lhd$ be an arbitrary linear order on $[1,n]$.
Then the elements $E_a$ for $a \in \ZZ^m$ form a basis in $\AA$ as a $\ZZ[v, v^{-1}]$-module, which satisfies the condition
\eqref{eq:bar-triangularity} with respect to the bar-involution on $\AA$, and the partial order on the index set $L = \ZZ^m$ given by
\eqref{eq:partial-order-on-Zm}.
Therefore, this basis gives rise to a $\ZZ[v, v^{-1}]$-basis $\BB = \BB(\tilde \XX, \tilde B; \lhd) = \{C_a \, : \, a \in \ZZ^m\}$ in $\AA$ uniquely determined by
\eqref{eq:b-bar-fixed} and \eqref{eq:b-triangularity}.
\end{theorem}

We refer to the basis $\BB = \BB(\tilde \XX, \tilde B; \lhd)$ in Theorem~\ref{th:Lusztig-lemma-ca} as the \emph{triangular basis} associated with an acyclic
quantum seed $(\tilde \XX, \tilde B)$, and a linear order $\lhd$ on $[1,n]$.
Theorem~\ref{th:Lusztig-lemma-ca} will be proved in Section~\ref{sec:proof-Lusztig-lemma-ca}.

\begin{remark}
Theorem~\ref{th:Lusztig-lemma-ca} allows the following generalization to the case of an arbitrary quantum cluster algebra $\AA$ (not necessarily possessing an acyclic quantum seed) and its arbitrary quantum seed $(\tilde \XX, \tilde B)$. 
Namely, we can still consider a \emph{lower bound} $\underline{\AA}(\tilde \XX, \tilde B)$, that is, 
a $\ZZ \PP$-subalgebra of $\AA$ generated by the elements $X_k, X'_k$ for $k \in [1,n]$ (cf. \eqref{eq:lower-qca}).
Clearly, $\underline{\AA}(\tilde \XX, \tilde B)$ is invariant under the bar-involution. 
As shown in \cite{bz-advances}, this subalgebra is equal to $\AA$ if and only if the quantum seed $(\tilde \XX, \tilde B)$ is acyclic.
Let $\underline{L} \subseteq L = \ZZ^m$ be the set of vectors $a \in \ZZ^m$ such that the full directed subgraph of $\Gamma(B)$ with the set of vertices 
$\{j \in [1,n] : a_j < 0\}$ is acyclic.  
It is not hard to deduce from the results in \cite{bz-advances} that the elements $E_a$ for $a \in \underline{L}$ form a basis in $\underline{\AA}(\tilde \XX, \tilde B)$ as a $\ZZ[v, v^{-1}]$-module.
Then Theorem~\ref{th:Lusztig-lemma-ca} as well as its proof in Section~\ref{sec:proof-Lusztig-lemma-ca} remain valid if one replaces $\AA$ with 
$\underline{\AA}(\tilde \XX, \tilde B)$, and $L = \ZZ^m$ with $\underline{L}$ (with the restricted partial order).  
\end{remark}

Varying an acyclic quantum seed $(\tilde \XX, \tilde B)$ and a linear order $\lhd$, we obtain  many different triangular bases in $\AA$.
However it turns out that many of these bases actually coincide with each other, and so provide us with a more ``canonical" choice.
To describe it, we recall the following characterization of acyclic seeds:
\begin{align}
\label{eq:acyclic-thru-linear-order}
&\text{A quantum seed $(\tilde \XX, \tilde B)$  is acyclic if and only if there exists a linear}\\
\nonumber
 &\text{order $\lhd$ on $[1,n]$ such that
$b_{ij} \leq 0$ for all $i, j \in [1,n]$ with $i \lhd j$}
\end{align}
(this follows at once from the well-known fact that every partial order on a finite set extends to a linear order).

Now we can state the main result of this paper.

\begin{theorem}
\label{th:canonical-lusztig-type-basis}
A triangular basis $\BB = \BB(\tilde \XX, \tilde B; \lhd)$ in a quantum cluster algebra $\AA$ does not depend on the choice of an acyclic quantum seed
$(\tilde \XX, \tilde B)$ and a linear order $\lhd$ on $[1,n]$ provided this linear order satisfies the condition in \eqref{eq:acyclic-thru-linear-order}.
\end{theorem}

We refer to the basis $\BB$ in Theorem~\ref{th:canonical-lusztig-type-basis} as the \emph{canonical triangular basis}.
Theorem~\ref{th:canonical-lusztig-type-basis} will be proved in Sections~\ref{sec:proof-main} - \ref{sec:proof-main-principal}.

By the construction, $\BB$ contains all elements $X^a$ for $a = (a_1, \dots, a_m) \in \ZZ^m$ such that $a_k \geq 0$ for $k \in [1,n]$.
We refer to these elements as (normalized) \emph{cluster monomials} associated with a quantum seed $(\tilde \XX, \tilde B)$.
As an immediate corollary of Theorem~\ref{th:canonical-lusztig-type-basis}, we obtain the following important property of~$\BB$.

\begin{corollary}
\label{cor:B-contains-acyclic-cluster-monomials}
The canonical triangular basis in $\AA$ contains all cluster monomials associated to acyclic quantum seeds of $\AA$.
\end{corollary}

\begin{remark}
F.~Qin (private communication) has informed us that the results of \cite{kimura-qin} imply the following: if $\AA$ has an acyclic quantum seed such that
the exchange matrix $B$ is skew-symmetic, and the directed graph $\Gamma(B)$ is \emph{bipartite} then $\BB$ contains \emph{all} cluster monomials associated to quantum seeds of ~$\AA$ (not only acyclic ones).
Furthermore, the results in \cite{GLS} imply that each acyclic quantum cluster algebra $\AA$ (with an appropriate choice of frozen variables)
can be realized as the localization with respect to frozen variables of the coordinate ring of a quantum unipotent cell (in an appropriate Kac-Moody group) associated with the square of a Coxeter element. 
Then the technique developed in \cite{GLS,kimura} implies that, under the same assumptions as above ($B$ skew-symmetric, and $\Gamma(B)$ bipartite), 
the canonical triangular basis  becomes identified with the dual canonical basis for the quantum unipotent cell (with the appropiate relation between the formal variable $q$ in \cite{GLS,kimura}, and the formal variable $v$ used in this paper).
\end{remark}


The rest of the paper is organized as follows.
In Section~\ref{sec:proof-Lusztig-lemma-ca} we prove Theorem~\ref{th:Lusztig-lemma-ca}.
Sections~\ref{sec:proof-main} - \ref{sec:proof-main-principal} are devoted to the proof of  Theorem~\ref{th:canonical-lusztig-type-basis}.

In Section~\ref{sec:proof-main} we reduce the proof of Theorem~\ref{th:canonical-lusztig-type-basis} to a special case
of \emph{principal quantization} described in Example~\ref{ex:principal-q-seed}.
This is done by a two-step construction that we find of independent interest.
First, we embed any quantum cluster algebra $\AA (B)$ with an acyclic exchange matrix $B$ into a ``double" algebra $\AA^{(2)}(B)$ obtained by adjoining some frozen variables
in an appropriate way; and then we show that the principal quantization $\AA_\bullet (B)$ can also be embedded in $\AA^{(2)}(B)$.
Note that the principal quantization appears (in some disguise) in \cite{fg} in the context of the ``quantum symplectic double" (we are grateful to A.~Goncharov for bringing
this connection to our attention).
It remains to be seen whether the above reduction procedure makes sense in the context of \cite{fg}.

The case of the principal quantization is treated in Section~\ref{sec:proof-main-principal}.
For the convenience of the reader we treat a special case $n=2$ separately in Section~\ref{sec:proof-main-principal-rank2}, to illustrate the general idea
in a somewhat less technical situation.

In Section~\ref{sec:coef-free-rank2-affine} we illustrate our results in the special case when $m=n=2$ and
$\tilde B = B =
\left(\!\!\begin{array}{cc}
0 & -2 \\
2 & 0  \\
\end{array}\!\!\right) \ .$
This cluster algebra and various bases in it were studied in detail in \cite{ding-xu,lampe}.
We make use of the calculations in \cite{ding-xu}, and we show (in Proposition~\ref{pr:can-triangular-basis-A11})
that in this special case our canonical triangular basis coincides with the dual canonical basis studied in \cite{lampe}, and also with one of the bases constructed in \cite{ding-xu}
with the help of the quantum Caldero-Chapoton characters
associated with indecomposable representations of the Kronecker quiver (see \cite{rupel}).

The concluding Section~\ref{sec:proof-Lusztig-lemma} contains the proof of Lusztig's Lemma (Theorem~\ref{th:Lusztig-lemma}).

\section{Proof of Theorem~\ref{th:Lusztig-lemma-ca}}
\label{sec:proof-Lusztig-lemma-ca}

In this section we prove Theorem~\ref{th:Lusztig-lemma-ca}.
We freely use the terminology and notation from Section~\ref{sec:intro}.
Recall that we work with the quantum cluster algebra $\AA = \AA(\tilde \XX, \tilde B)$ associated with an acyclic quantum seed.
According to \eqref{eq:lower-qca}, $\AA$ is a $\ZZ \PP$-algebra
generated by the elements $X_k, X'_k$ for $k \in [1,n]$.
We start by collecting the quasi-commutation relations among these generators.
The relations between the generators $X_k$ are given by \eqref{eq:Xi-q-com}.
The rest are as follows (recall that the vectors $e'_k$ were introduced in \eqref{eq:e-prime-k}).

\begin{lemma}
\label{lem:commut-relations}
\begin{enumerate}
\item For $i \in [1,m]$ and $k \in [1,n]$ with $i \neq k$ we have
\begin{equation}
\label{eq:Xi-Xk-prime-com}
X_i X'_k = v^{2 \Lambda(e_i,e'_k)} X'_k X_i \, .
\end{equation}
\item For $k \in [1,n]$  we have
\begin{equation}
\label{eq:Xk-Xk-prime-com}
v^{-\Lambda(e'_k,e_k)}X'_k X_k  - v^{\Lambda(e'_k,e_k)}X_k X'_k = (v^{-d_k} - v^{d_k}) X^{[-b_k]_+} \, .
\end{equation}
\item For $j, k \in [1,n]$ with $j \neq k$ we have
\begin{equation}
\label{eq:Xj-prime-Xk-prime-com}
v^{-\Lambda(e'_j,e'_k)}X'_j X'_k  - v^{\Lambda(e'_j,e'_k)}X'_k X'_j = (v^{-d_j b_{jk}} - v^{d_j b_{jk}})
X^{-e_j - e_k + [-\varepsilon_{jk} b_j]_+ + [\varepsilon_{jk} b_k]_+} \, ,
\end{equation}
where $\varepsilon_{jk}$ is the sign of the matrix entry $b_{jk}$.
\end{enumerate}
\end{lemma}

\begin{proof}
Remembering the notation \eqref{eq:e-prime-k}, we can rewrite \eqref{eq:q-exchange-rel-xx} as
\begin{equation}
\label{eq:X-prime-k}
X'_k = X^{e'_k} \ + \ X^{e'_k -b_{k}} \ .
\end{equation}
Now the equalities \eqref{eq:Xi-Xk-prime-com} - \eqref{eq:Xj-prime-Xk-prime-com} follow easily from
\eqref{eq:multiplication-quantum-torus} and \eqref{eq:orthogonality}
(note that \eqref{eq:Xj-prime-Xk-prime-com} was essentially proved in \cite[Proposition~7.1]{bz-advances}).
\end{proof}

Consider an increasing filtration $\{0\} = F_{-1} \subset F_0 \subset F_1 \subset \cdots$ on $\AA$ defined as follows:
\begin{align}
\nonumber
&\text{$F_r$ is the $\ZZ\PP$-linear span of (noncommutative) monomials $M$}\\
\label{eq:filtration}
&\text{in the elements $X_k$ and $X'_k$ for $k \in [1,n]$
such that the degree}\\
\nonumber
&\text{of $M$ with respect to the elements $X'_k$ does not exceed $r$.}
\end{align}

By the definition we have $F_r \cdot F_s \subseteq F_{r+s}$ for $r, s \geq 0$.
Thus, this filtration gives rise to an \emph{associated graded $\ZZ \PP$-algebra} $\widehat \AA$ given by
\begin{equation}
\label{eq:assoc-graded}
\text{$\widehat \AA = \bigoplus_{r \geq 0} \widehat {\AA}_r$, where $\widehat {\AA}_r = F_r/F_{r-1}$.}
\end{equation}
For an element $E \in F_r - F_{r-1}$ we denote by $\widehat E$ the (non-zero) image of $E$ in $\widehat {\AA}_r$.

The relations in Lemma~\ref{lem:commut-relations} imply that \eqref{eq:monomial-basis-lower-qca} can be sharpened as follows:
\begin{align}
\label{eq:monomial-basis-Fr}
&\text{For every $r \geq 0$, the elements $E_a^\circ$ with $r(a) \leq r$}\\
\nonumber
&\text{form a basis in $F_r$ as a $\ZZ[v, v^{-1}]$-module}
\end{align}
(see \eqref{eq:r-of-a}).
It follows that
\begin{align}
\label{eq:monomial-basis-Fr-graded}
&\text{For every $r \geq 0$, the elements $\widehat {E_a^\circ}$ with $r(a) = r$}\\
\nonumber
&\text{form a basis in $\widehat {\AA}_r$ as a $\ZZ[v, v^{-1}]$-module.}
\end{align}
Clearly, \eqref{eq:monomial-basis-Fr} and \eqref{eq:monomial-basis-Fr-graded} remain true if we replace
$E_a^\circ$ with $E_a = v^{\nu(a)} E_a^\circ$.

Now we turn our attention to the bar-involution on $\AA$.
In view of \eqref{eq:filtration} and  \eqref{eq:bar-properties},
we have $\overline {F_r} = F_r$ for all~$r$.
Therefore, this involution gives rise to an involution on each $\widehat {\AA}_r$, which we will denote in the same way.
In this notation, the desired condition \eqref{eq:bar-triangularity} for the partial order on the index set $L = \ZZ^m$ given by
\eqref{eq:partial-order-on-Zm}, simply means that every element $\widehat {E_a}$ is invariant under the bar-involution.

Fix some $a \in \ZZ^m$.
To prove the desired equality $\overline {\widehat {E_a}} = \widehat {E_a}$, we can view $\widehat {E_a}$
as an element of a new quantum torus $\TT(a)$ generated by the following elements:
\begin{itemize}
\item
the frozen variables $X_i$ for $i \in [n+1,m]$;
\item
the elements $X_j$ for $j \in [1,n]$ such that $a_j > 0$; and
\item
the elements $\widehat {X'_k}$ for $k \in [1,n]$ such that $a_k < 0$.
\end{itemize}
The commutation relations among the generators of the first two kinds are given by
\eqref{eq:Xi-q-com}, and, according to \eqref{eq:Xi-Xk-prime-com} and \eqref{eq:Xj-prime-Xk-prime-com}, those involving one or both generators
$\widehat {X'_k}$ are of the form
\begin{equation}
\label{eq:Xk-prime-com-graded}
X_i \widehat {X'_k} = v^{2 \Lambda(e_i,e'_k)} \widehat {X'_k} X_i, \quad
\widehat {X'_j}  \widehat {X'_k}  =  v^{2\Lambda(e'_j,e'_k)} \widehat {X'_k} \widehat {X'_j} \ ;
\end{equation}
thus, they are of the same form \eqref{eq:Xi-q-com}, with the following modification: the argument of the form $\Lambda$
(that appears in the exponent of~$v$), corresponding to each $\widehat {X'_k}$ is equal to $e'_k$.
We see that the condition that $\widehat {E_a}$ is invariant under the bar-involution, is in fact equivalent
to the definition of $E_a$ via \eqref{eq:E-a} and \eqref{eq:nu-of-a}.
This completes the proof of Theorem~\ref{th:Lusztig-lemma-ca}.

The tools developed in the course of the above proof  allow us to show the following.

\begin{proposition}
\label{pr:B-independent-under-transposition}
For a given acyclic quantum seed $(\tilde \XX, \tilde B)$, the elements $E_a$ associated with a linear order $\lhd$ on $[1,n]$, do not change under the
following modification of $\lhd$: interchanging two adjacent indices $j$ and $k$ such that $b_{jk} = 0$.
Therefore, the triangular basis $\BB(\tilde \XX, \tilde B; \lhd)$ also does not change under such a modification.
\end{proposition}

\begin{proof}
If $b_{jk} = 0$ then the right hand side of \eqref{eq:Xj-prime-Xk-prime-com} becomes equal to~$0$, hence the elements $X'_j$ and $X'_k$
quasi-commute (that is, commute, up to a multiple which is a power of~$v$).
Remembering \eqref{eq:standard-monomial-non-normalized}, we see that
interchanging $j$ and $k$ results in multiplying each $E_a^\circ$ by a power of~$v$.
But the element $E_a$ can be characterized as an element of the form $v^s E_a^\circ$ such that
$\widehat {E_a}$ is invariant under the bar-involution.
Thus, $E_a$ remains unchanged under the modification in question.
\end{proof}

\section{Proof of Theorem~\ref{th:canonical-lusztig-type-basis}-I. Reduction to principal case}
\label{sec:proof-main}

In this and the next two sections we prove Theorem~\ref{th:canonical-lusztig-type-basis}, that is, that the
triangular  basis $\BB = \BB(\tilde \XX, \tilde B; \lhd)$ in a quantum cluster algebra $\AA$ does not depend on the choice of an acyclic quantum seed
$(\tilde \XX, \tilde B)$ and a linear order $\lhd$ on $[1,n]$ as long as  $\lhd$ satisfies the condition in \eqref{eq:acyclic-thru-linear-order}.

First we fix an acyclic quantum seed $(\tilde \XX, \tilde B)$ and show that $\BB$ does not depend
on the choice of a linear order $\lhd$ satisfying \eqref{eq:acyclic-thru-linear-order}.
Indeed, since the directed graph $\Gamma(B)$ has no oriented cycles, it gives rise to the partial order on $[1,n]$
defined as follows: $i$ is less than $j$ if there is an oriented path from $i$ to $j$ in $\Gamma(B)$.
Now the condition \eqref{eq:acyclic-thru-linear-order} means that $\lhd$ is a linear extension of this partial order.
It is well-known (and easy to prove) that every two linear extensions of a given finite partial order can be obtained from each other by a sequence of operations of the following kind:
interchanging two adjacent elements provided they are incomparable with respect to the partial order in question.
Thus, the desired independence of $\BB$ on $\lhd$ follows from Proposition~\ref{pr:B-independent-under-transposition}.

\smallskip

It remains to show that $\BB$ is independent of the choice of an acyclic quantum seed $(\tilde \XX, \tilde B)$.
As shown in \cite[Corollary~4]{CK}, every two acyclic quantum seeds
can be obtained from each other by a finite number of \emph{shape-preserving} mutations (recall that a mutation $\mu_k$
is shape-preserving if $k \in [1,n]$ is a source or a sink of the directed graph $\Gamma(\tilde B)$). 
\footnote{In \cite{CK} the exchange matrix~$B$ was assumed to be skew-symmetric; the general case can be deduced from this one by the standard argument involving
\emph{folding}.}
Therefore, to finish the proof of Theorem~\ref{th:canonical-lusztig-type-basis}, it is enough to show that the basis
$\BB(\tilde \XX, \tilde B; \lhd)$ with the order $\lhd$ satisfying \eqref{eq:acyclic-thru-linear-order}, does not change under a shape-preserving mutation.
Without loss of generality, we assume for the rest of this section that the order $\lhd$ is the usual one: $1 \lhd 2 \lhd \cdots \lhd n$.
Thus, we have
\begin{equation}
\label{eq:signs-B-natural-order}
\text{$b_{ij} \leq 0$ for $1 \leq i < j \leq n$.}
\end{equation}
For a vector $a = (a_1, \dots, a_m) \in \ZZ^m$, we abbreviate
\begin{equation}
\label{eq:truncated-vectors}
a^{> n} = \sum_{n < i \leq m} a_i e_i, \quad a^{\leq n} = \sum_{1 \leq i \leq n} a_i e_i \ .
\end{equation}
Note that the truncation operations $a \mapsto a^{> n}$ and $a \mapsto a^{\leq n}$ commute with the operation
$a \mapsto [a]_+$ (replacing each component $a_i$ with $[a_i]_+ = \max(a_i,0)$), hence we have well defined
operations $a \mapsto [a]_+^{> n}$ and  $a \mapsto [a]_+^{\leq n}$.
In this notation, the definition \eqref{eq:standard-monomial-non-normalized} takes the form
\begin{equation}
\label{eq:standard-monomial-standard-order}
E_a^\circ = X^{a^{>n} + [a]_+^{\leq n}} \cdot
(X'_1)^{[- a_1]_+} (X'_2)^{[- a_2]_+} \cdots (X'_n)^{[- a_n]_+}\ .
\end{equation}

Again without loss of generality, as a shape-preserving mutation, we take $\mu_n$, the mutation at the last index~$n$.
The mutated quantum seed $\mu_n(\tilde \XX, \tilde B) =  (\tilde {\XX'}, \tilde {B'})$ is defined as follows.
The new extended cluster $\tilde {\XX'}$ is obtained from $\tilde \XX = \{X_1, \dots, X_m\}$ by replacing $X_n$ with $X'_n$ given by
\eqref{eq:q-exchange-rel-xx} or \eqref{eq:X-prime-k}  (for $k=n$), where the vector $e'_n$ is given by \eqref{eq:e-prime-k}.
In view of \eqref{eq:signs-B-natural-order}, we have
\begin{equation}
\label{eq:e-prime-n-special}
e'_n = - e_n + b_n^{> n} \, ,
\end{equation}
The set $\tilde {\XX'}$ generates a based quantum torus $\TT' \subset \FF$ defined in the same way as the quantum torus $\TT$ above, but with the quasi-commutation relations
governed by the new skew-symmetric form $\Lambda'$ on $\ZZ^m$ given as follows: $\Lambda'(e_i, e_j) = \Lambda(e_i, e_j)$ for $i, j \neq n$, and
$\Lambda'(e_n, e_j) = \Lambda (e'_n, e_j)$.

Recall that the new extended exchange matrix $\tilde B'$ is obtained from $\tilde B$ by the \emph{matrix mutation} in \eqref{eq:matrix-mutation}
with $k = n$.
Note that the sign condition $b_{ij} \leq 0$ for $1 \leq i < j \leq n$ implies that the top $n \times n$ part $B'$ of $\tilde {B'}$ is obtained from the top part $B$ of
$\tilde B$ by simply changing the sign of all entries in the $n$th row and column.
It follows that the condition \eqref{eq:acyclic-thru-linear-order} for the new seed will be obeyed by choosing the linear order $\lhd$ on $[1,n]$ as follows:
\begin{equation}
\label{eq:shifted-linear-order}
n \lhd 1 \lhd 2 \lhd \cdots \lhd (n-1) \ .
\end{equation}

For every $a \in \ZZ^m$, we define the element $(E'_a)^\circ \in \AA$ is the same way as in \eqref{eq:standard-monomial-non-normalized} but with respect to the
acyclic quantum seed $(\tilde {\XX'}, \tilde {B'})$ and the linear order $\lhd$ given by \eqref{eq:shifted-linear-order}.
To describe $(E'_a)^\circ$ explicitly, we need some more notation.
Define the vectors $e''_1, \dots, e''_n$ by setting
\begin{equation}
\label{eq:e-double-prime-k}
e''_k = - e_k + [b'_k]_+ \in \ZZ^m \, ,
\end{equation}
where $b'_k$ is the $k$th column of $\tilde B'$.
Then the mutated cluster variables from the quantum seed $(\tilde {\XX'}, \tilde {B'})$
are given by the following counterpart
of \eqref{eq:X-prime-k}:
\begin{equation}
\label{eq:X-double-prime-k}
X''_k = (X')^{e''_k} \ + \ (X')^{e''_k -b'_{k}} \ .
\end{equation}
Note that we have $X''_n = X_n$ since the mutation $\mu_n$ is an involution; however it is sometimes convenient
to keep the notation $X''_n$.

Using this notation, we can obtain $(E'_a)^\circ$ from the element $E_a^\circ$ given by \eqref{eq:standard-monomial-standard-order} by the following modifications:
\begin{itemize}


\item the factor $X^{a^{>n} + [a]_+^{\leq n}}$ is replaced with its counterpart
$$(X')^{a^{>n} + [a]_+^{\leq n}}$$
in the quantum torus $\TT'$ (that is, in constructing this element,
$X_n$ is replaced by $X'_n$).

\item the remaining product $(X'_1)^{[- a_1]_+} (X'_2)^{[- a_2]_+} \cdots (X'_n)^{[- a_n]_+}$ is replaced by
$$(X''_n)^{[- a_n]_+} (X''_1)^{[- a_1]_+} (X''_2)^{[- a_2]_+} \cdots (X''_{n-1})^{[- a_{n-1}]_+} \, .$$
\end{itemize}

As for the normalized basis element $E'_a$, it is obtained from $(E'_a)^\circ$ by a suitable modification of \eqref{eq:E-a} and
\eqref{eq:nu-of-a}.
Namely, we have $E'_a = v^{\nu'(a)} (E'_a)^\circ$, with the exponent $\nu'(a) \in \ZZ$ determined from the condition that
$v^{\nu'(a)} LT((E'_a)^\circ)$ is invariant under the bar-involution, where
the leading term $LT((E'_a)^\circ)$ is obtained from
$(E_a)^\circ$ by replacing each $X''_k$ with $(X')^{e''_k}$.

In particular, we have:
\begin{eqnarray}
\label{eq:E-prime-units-easy}
& E'_{\pm e_i} = E_{\pm e_i} = X_i^{\pm 1} \,\, (n < i \leq m), \quad
& E'_{e_k}  = E_{e_k} = X_k \,\, (1 \leq k < n) \, ,\\
\nonumber
& E'_{e_n} = E_{-e_n} = X'_n, \quad & E'_{-e_n}  = E_{e_n} = X_n \ .
\end{eqnarray}
Also for $1 \leq k < n$, we have $E'_{-e_k} = X''_k$.
Let us compute the expansion of this element in the basis  $\{E_a\}$.

Recall that, for $r \geq s \geq 0$, the Gaussian binomial coefficient $\begin{bmatrix} r \\ s \\ \end{bmatrix}_t$ is an integer polynomial in a variable~$t$ given by
\begin{equation}
\label{eq:gaussian-binom}
\begin{bmatrix} r \\ s \\ \end{bmatrix}_t = \frac{(t^r - 1)(t^{r-1} - 1) \cdots (t^{r-s+1} - 1)}{(t^s - 1)(t^{s-1} - 1) \cdots (t-1)} \, .
\end{equation}
The corresponding binomial formula is
\begin{equation}
\label{eq:gaussian-binom-identity}
\prod_{p=1}^r (1 + u^{2p-1} X) = \sum_{s=0}^r u^{s^2} \begin{bmatrix} r \\ s \\ \end{bmatrix}_{u^2} X^s \ .
\end{equation}

Now, for $1 \leq k < n$, we define a vector $\varphi(-e_k) \in \ZZ^m$ as follows (the notation will be explained in a moment):
\begin{equation}
\label{eq:Xk-double-prime-leading-term}
\varphi(-e_k) = - e_k - b_{nk} e_n + [-b'_{k}]_+^{>n} - [-b_k]_+^{>n} \ .
\end{equation}
Then we have
\begin{equation}
\label{eq:Xk-double-prime-expansion}
 E'_{-e_k} = X''_k = E_{\varphi(-e_k)} - \sum_{s=1}^{b_{nk}} v^{s^2 d_n} \begin{bmatrix} b_{nk} \\ s \\ \end{bmatrix}_{v^{2d_n}}
 E_{e''_k - s b_n} \, .
 \end{equation}
 This calculation was essentially done in \cite[Proof of Lemma~5.8]{bz-advances}.
 Note that if $b_{nk} = b_{kn} = 0$ then the last sum in \eqref{eq:Xk-double-prime-expansion} disappears, and $\varphi(-e_k) = -e_k$, so in this case
 $X''_k = E'_{-e_k} = E_{-e_k} = X'_k$.
 And if $b_{kn} < 0$ (hence $b_{nk} > 0$) then, for every $j \in [1,n]$ and $s \in [1,b_{nk}]$, the $j$th component of the vector
 $e''_k - s b_n$ is nonnegative, therefore we have $E_{e''_k - s b_n} = X^{e''_k - s b_n}$.

 Let $\AA_+ = \AA_+(\tilde \XX, \tilde B)$ denote the $\ZZ[v]$-linear span of the basis $\{E_a\}$.
 We refer to $\AA_+$ as the \emph{crystal lattice} associated with a given acyclic quantum seed.
 Our proof of Theorem~\ref{th:canonical-lusztig-type-basis} is based on the following fact.

 \begin{theorem}
 \label{th:E-primes-crystal-lattice}
 All elements $E'_a$ belong to $\AA_+ - v \AA_+$. More precisely, there is a bijection $\varphi: \ZZ^m \to \ZZ^m$ such that, for every $a \in \ZZ^m$, we have
 \begin{equation}
 \label{eq:E'-E-crystal}
 E'_a - E_{\varphi(a)} \in v \AA_+ \ .
 \end{equation}
\end{theorem}

Before proving Theorem~\ref{th:E-primes-crystal-lattice}, we deduce the following corollary.

\begin{corollary}
\label{cor:B-prime-thru-B}
The triangular basis $\BB' = \BB(\tilde {\XX'}, \tilde {B'})$ coincides with
$\BB = \BB(\tilde {\XX}, \tilde {B})$.
More precisely, for any $a \in \ZZ^m$, an element $C'_a$ of $\BB'$ is equal to $C_{\varphi(a)}$.
\end{corollary}

\begin{proof}
Theorem~\ref{th:E-primes-crystal-lattice} implies that the crystal lattice $\AA'_+ = \AA_+ (\tilde {\XX'}, \tilde {B'})$
is contained in $\AA_+$.
Thus, the condition \eqref{eq:b-triangularity} for the basis $\BB'$ implies that
$$C'_a  \in E'_a + v \AA'_+ \subseteq E_{\varphi(a)} + v \AA_+ \ .$$
In view of Theorem~\ref{th:Lusztig-lemma-ca}, we have $C'_a = C_{\varphi(a)}$, as claimed.
\end{proof}

In view of \eqref{eq:Xk-double-prime-expansion}, the condition \eqref{eq:E'-E-crystal} holds for $a = - e_k$ with $k < n$, and $\varphi(- e_k)$ given by  \eqref{eq:Xk-double-prime-leading-term}.
And in view of \eqref{eq:E-prime-units-easy}, it also holds for the rest of the vectors $\pm e_i$, with  $\varphi (\pm e_i)$ given as follows:
\begin{eqnarray}
\label{eq:varphi-easy-units}
& \varphi(\pm e_i) = \pm e_i \,\, (n < i \leq m), \quad
& \varphi(e_k) = e_k \,\, (k < n) \, ,\\
\nonumber
& \varphi(e_n) = -e_n, \quad & \varphi(-e_n)  = e_n \ .
\end{eqnarray}


For every $a \in \ZZ^m$, we write the element $E'_a$ in terms of the basis $\{E_{a'}\}$ as follows:
\begin{equation}
\label{eq:E-prime-E=expansion}
E'_a = \sum_{a' \in \ZZ^m} c^{a'}_a  E_{a'} \ ,
\end{equation}
where all coefficients $c^{a'}_a$ belong to $\ZZ[v^{\pm 1}]$, and, for a given~$a$, all but finitely many of them
are equal to~$0$.
It is easy to see that to prove Theorem~\ref{th:E-primes-crystal-lattice}, it suffices to show the following:
\begin{align}
\label{eq:c-conditions}
&\text{For every $a \in \ZZ^m$, exactly one of the coefficients}\\
\nonumber
&\text{$c^{a'}_a$ is equal to $1$, and the rest
belong to $v \ZZ[v]$}.
\end{align}

The proof of \eqref{eq:c-conditions} will occupy the rest of this section as well as the next two sections.
We start with the following lemma.

\begin{lemma}
\label{lem:shift-by-coeff}
We have $c^{a'+ a_\circ}_{a + a_\circ} = c^{a'}_a$ for all $a, a', a_\circ \in \ZZ^m$ such that $a_\circ^{\leq n} = 0$.
\end{lemma}

\begin{proof}
The key observation is that the element $(X')^{a_\circ}$ is equal to $X^{a_\circ}$, and that it quasi-commutes (that is, commutes up to an integer
power of $v$ as a multiple) with all elements $E'_a$ and $E_{a'}$; furthermore, $X^{a_\circ}$ has the same quasi-commutation multiple with $E'_a$ and
$LT((E'_a)^\circ)$, as well as with $E_{a'}$ and $LT(E_{a'}^\circ)$.
Then the definitions readily imply that
\begin{equation}
\label{eq:Xacirc-times-Ea}
X^{a_\circ} E'_a = v^{\nu'(a_\circ; a)} E'_{a + a_\circ}, \quad
X^{a_\circ} E_{a'} = v^{\nu(a_\circ; a')} E_{a' + a_\circ} \, ,
\end{equation}
where the integer exponents $\nu'(a_\circ; a)$ and $\nu(a_\circ; a')$ are determined from the quasi-commutation relations
\begin{equation}
\label{eq:Xacirc-times-Ea-q-commute}
X^{a_\circ} E'_a = v^{2 \nu'(a_\circ; a)}  E'_a X^{a_\circ}, \quad
X^{a_\circ} E_{a'} = v^{2 \nu(a_\circ; a')}  E_{a'} X^{a_\circ}  \, .
\end{equation}
Now the fact that $X^{a_\circ}$ quasi-commutes with all the terms in the equality \eqref{eq:E-prime-E=expansion} implies that
the corresponding quasi-commutation multiples are all the same, that is, $\nu'(a_\circ; a) = \nu(a_\circ; a')$ whenever
$c^{a'}_a \neq 0$.
Thus, the desired equality $c^{a'+ a_\circ}_{a + a_\circ} = c^{a'}_a$ follows from \eqref{eq:Xacirc-times-Ea}.
\end{proof}

In view of Lemma~\ref{lem:shift-by-coeff}, it is enough to prove \eqref{eq:c-conditions} for $a \in \ZZ^n$.
This will be done in two steps.

The first step (occupying the rest of this section) is a reduction to a special choice of a quantum seed $(\tilde \XX^\bullet, \tilde B^\bullet)$ (for a given exchange matrix~$B$) namely, to
the \emph{principal quantum seed} introduced in Example~\ref{ex:principal-q-seed}.
To do this, we embed the quantum cluster algebra $\AA(\tilde \XX, \tilde B)$ into a bigger one associated with the ``double" quantum seed
$(\tilde \XX^{(2)}, \tilde B^{(2)})$ defined by the following data:
\begin{itemize}
\item
$\tilde B^{(2)}$ is an integer $2m \! \times \! n$ matrix with the top $m \! \times \! n$ block $\tilde B$, and the bottom $m \! \times \! n$
block a zero matrix.
 \item
A skew-symmetric bilinear form $\Lambda^{(2)}: \ZZ^{2m} \times \ZZ^{2m} \to \ZZ$  is given by
\begin{equation}
\label{eq:Lambda-double}
\Lambda^{(2)}((e,e'), (f,f'))= \Lambda(e,f) - \Lambda(e',f') \quad (e, e', f, f' \in \ZZ^m) \, .
\end{equation}
\end{itemize}
Note that $\tilde B^{(2)}$ and $\Lambda^{(2)}$ satisfy the compatibility condition \eqref{eq:orthogonality} with the same choice
of positive integers $d_1, \dots, d_n$.

The corresponding \emph{double based quantum torus}
$\TT^{(2)} = \TT(\Lambda^{(2)})$ is the $\ZZ[v, v^{-1}]$-algebra
with a distinguished $\ZZ[v, v^{-1}]$-basis $\{X^{(e,e')}: e, e' \in \ZZ^m\}$ and the
multiplication given by
\begin{equation}
\label{eq:multiplication-quantum-torus-double}
X^{(e,e')} X^{(f,f')} = v^{\Lambda^{(2)}((e,e'), (f,f'))} X^{e+e', f + f'} \quad (e, e', f, f' \in \ZZ^m) \, .
\end{equation}
We identify the based quantum torus $\TT = \TT(\Lambda)$ with a subalgebra in $\TT^{(2)}$ by setting
$X^e  = X^{(e,0)}$ for $e \in \ZZ^m$.
This way the quantum cluster algebra $\AA(\tilde \XX, \tilde B)$ is identified with a subalgebra of $\AA(\tilde \XX^{(2)}, \tilde B^{(2)})$.

For $a^{(2)} \in \ZZ^{2m}$, we denote by $E^{(2)}_{a^{(2)}}$ the element of $\AA(\tilde \XX^{(2)}, \tilde B^{(2)})$ constructed as in \eqref{eq:standard-monomial-non-normalized}
and \eqref{eq:E-a}.
In particular, we have
\begin{equation}
\label{eq:Ea-double}
E_a = E^{(2)}_{(a,0)}
\end{equation}
for each $a \in \ZZ^m$.

It turns out that the principal quantum cluster algebra $\AA(\tilde \XX^\bullet, \tilde B^\bullet)$
can also be realized as a subalgebra of $\AA(\tilde \XX^{(2)}, \tilde B^{(2)})$.
To avoid notational confusion, we will use the symbol $\bullet$ for the objects related to this subalgebra.
First we introduce the based quantum torus $\TT^\bullet$ as the $\ZZ[v, v^{-1}]$-subalgebra
of $\TT^{(2)}$ generated by the elements $X^\bullet_1, \dots, X^\bullet_{2n}$ given by
\begin{equation}
\label{eq:X-bullet-i}
X^\bullet_j = X^{(e_j,e_j)}, \quad X^\bullet_{n+j} = X^{(b_j^{>n}, - b_j^{\leq n})}
\quad (j \in [1,n]) \, ,
\end{equation}
where we use the notation \eqref{eq:truncated-vectors} for truncated vectors, and as before, $b_j \in \ZZ^m$ is the $j$th column of $\tilde B$.

\begin{lemma}
\label{lem:principal-within-double}
The based quantum torus $\TT^\bullet$ is equal to $\TT(\Lambda^\bullet)$, where $\Lambda^\bullet$ is the skew-symmetric bilinear form
on $\ZZ^{2n}$ with the matrix \eqref{eq:lambda-principal}.
\end{lemma}

\begin{proof}
First we show that the $2n$ vectors
$$(e_1,e_1), \dots, (e_n, e_n), (b_1^{>n}, - b_1^{\leq n}), \dots, (b_n^{>n}, - b_n^{\leq n})$$
in $\ZZ^{2m}$ are linearly independent, implying that $X^\bullet_1, \dots, X^\bullet_{2n}$ are algebraically independent.
Indeed, for every $j \in [1,n]$, we have
$$(b_j^{>n}, - b_j^{\leq n}) + \sum_{i=1}^n b_{ij} (e_i,e_i) = (b_j, 0) \, .$$
Thus, it is enough to show that the vectors $(e_1,e_1), \dots, (e_n, e_n), (b_1, 0), \dots, (b_n, 0)$ are linearly independent.
This follows by noticing that the first $n$ of these vectors have linearly independent second components, and then that the last $n$ vectors
are linearly independent in view of \eqref{eq:orthogonality}.

It remains to show that the elements $X^\bullet_1, \dots, X^\bullet_{2n}$ obey the commutation relations in \eqref{eq:Xi-q-com-principal}, and also that
the first $n$ of them commute with each other.
The latter statement follows from \eqref{eq:multiplication-quantum-torus-double} and \eqref{eq:Lambda-double}.
To show that, for $i, j \in [1,n]$,  the elements $X^\bullet_{n+i}$ and $X^\bullet_{n+j}$ satisfy the second equality in \eqref{eq:Xi-q-com-principal},
it is enough to observe that
\begin{align*}
\Lambda^{(2)}((b_i^{>n}, - b_i^{\leq n}), (b_j^{>n}, - b_j^{\leq n})) &=& &\Lambda(b_i^{>n}, b_j^{>n}) - \Lambda(b_i - b_i^{>n}, b_j - b_j^{>n})\\
&=& &- \Lambda(b_i, b_j) = d_j b_{ji} \, ,
\end{align*}
where we used \eqref{eq:orthogonality} and \eqref{eq:skew-symm}.
The first equality in \eqref{eq:Xi-q-com-principal} is checked in a similar way.
\end{proof}

In view of Lemma~\ref{lem:principal-within-double}, the assignment \eqref{eq:X-bullet-i} allows us to identify the principal quantum cluster algebra $\AA(\tilde \XX^\bullet, \tilde B^\bullet)$
with a subalgebra of $\AA(\tilde \XX^{(2)}, \tilde B^{(2)})$.
In accordance with the above notational convention, we denote the analogues of elements $E_a$ and $E'_a$ in
in $\AA(\tilde \XX^\bullet, \tilde B^\bullet)$ as  $E^\bullet_{a^\bullet}$ and ${E^\bullet_{a^\bullet}}'$, where $a^\bullet$ runs over $\ZZ^{2n}$.

\begin{lemma}
\label{lem:E-principal-double}
There exist injective maps $\psi: \ZZ^{2n} \to \ZZ^{2m}$ and $\psi': \ZZ^{2n} \to \ZZ^{2m}$ such that, for every $a \in \ZZ^{2n}$, we have
$E^\bullet_{a} = E^{(2)}_{\psi(a)}, \, {E'}^\bullet_{a} = {E'}^{(2)}_{\psi'(a)}$.
Specifically these maps are uniquely determined by the following properties:
\begin{itemize}
\item Both $\psi$ and $\psi'$ restrict as linear maps to the sublattice $\{0\} \times \ZZ^n \subset \ZZ^{2n}$ generated by $e_{n+1}, \dots, e_{2n}$ and to each of the $2^n$
coordinate orthants in the sublattice $\ZZ^n \times \{0\}$ generated by $e_1, \dots, e_n$. Moreover, if $a, a_\circ \in \ZZ^{2n}$, and $a_\circ^{\leq n} = 0$ then
$$\psi(a + a_\circ) = \psi(a) + \psi(a_\circ), \quad \psi'(a + a_\circ) = \psi'(a) + \psi'(a_\circ) \ .$$
\item We have
\begin{align*}
& \psi(e_{n+k}) = \psi'(e_{n+k}) = (b_k^{>n}, - b_k^{\leq n}) \quad (k \in [1,n])\ , \\
& \psi(e_k) = (e_k, e_k), \quad \psi(-e_k) = (- e_k - [-b_k]_+^{> n}, - e_k + [-b_k]_+^{\leq n}) \quad (k \in [1,n]) \ , \\
& \psi'(e_k) = (e_k, e_k) \quad (k \in [1,n-1]), \quad \psi'(e_n) = (e_n - [-b_n]_+^{> n}, - e_n -b_n^{\leq n}) \ , \\
& \psi'(-e_k) = (-e_k - [-b'_k]_+^{> n} - b_{nk} \cdot [-b_n]_+^{> n}, -e_k + [-b_k]_+^{\leq n} - b_{nk} \cdot b_n^{\leq n} - b_{nk} e_n) \\
& (k \in [1,n-1]),  \quad \psi'(-e_n) = (-e_n, e_n) \, .
\end{align*}
\end{itemize}
In particular, for every $a \in \ZZ^{2n}$, we have $\psi(a)^{\leq n} = \psi'(a)^{\leq n} = a^{\leq n}$.
\end{lemma}

\begin{proof}
This is proved by a direct calculation that compares the expressions \eqref{eq:q-exchange-rel-xx} and  \eqref{eq:Xk-double-prime-expansion}
evaluated in the principal algebra  $\AA(\tilde \XX^\bullet, \tilde B^\bullet)$ with their counterparts in the double algebra $\AA(\tilde \XX^{(2)}, \tilde B^{(2)})$.
\end{proof}

We are finally in a position to show the desired reduction to the principal case, namely that the validity of
\eqref{eq:c-conditions} in $\AA(\tilde \XX^\bullet, \tilde B^\bullet)$ implies that in  $\AA(\tilde \XX, \tilde B)$.
Let $a \in \ZZ^n$.
In view of Lemma~\ref{lem:shift-by-coeff}, the coefficients $c_a^{a'}$ in the expansion of $E'_a$ in the basis $\{E_{a'}\}$ of
$\AA(\tilde \XX, \tilde B)$ coincide with the corresponding coefficients in the expansion of ${E'}^{(2)}_{\psi'(a)}$
in $\AA(\tilde \XX^{(2)}, \tilde B^{(2)})$.
Since ${E'}^{(2)}_{\psi'(a)} = {E'}^\bullet_{a}$, its expansion in $\AA(\tilde \XX^{(2)}, \tilde B^{(2)})$ coincides with the corresponding expansion
in $\AA(\tilde \XX^\bullet, \tilde B^\bullet)$.
Therefore, if one of the coefficients in the latter expansion is equal to $1$, while the rest of the coefficients belong to $v \ZZ[v]$, then
the same is true for the expansion of $E'_a$ in $\AA(\tilde \XX, \tilde B)$, as desired.

\section{Proof of Theorem~\ref{th:canonical-lusztig-type-basis}-II. Principal quantization in rank~$2$}
\label{sec:proof-main-principal-rank2}

In this and the next section we finish the proof of Theorem~\ref{th:canonical-lusztig-type-basis} by showing that Theorem~\ref{th:E-primes-crystal-lattice}
and the condition \eqref{eq:c-conditions} hold in the quantum cluster algebra $\AA(\tilde \XX^\bullet, \tilde B^\bullet)$ with principal coefficients.
To make the argument more clear, we first explain it in the case $n = 2$, avoiding heavy notation.

We fix positive integers $b$ and $c$, and let the initial extended exchange matrix $\tilde B$ be given by
\begin{equation}
\label{eq:tilde-B-principal-rank2}
\tilde B =
\left(\!\!\begin{array}{cc}
0 & -b \\
c & 0 \\
1 & 0 \\
0 & 1 \\
\end{array}\!\!\right)
\end{equation}
(the signs of entries in the first two rows are chosen in accordance with \eqref{eq:signs-B-natural-order}).
Then we can take
\begin{equation}
\label{eq:d1-d2-rank2}
d_1 = c, \quad d_2 = b \, ,
\end{equation}
and the matrix $\Lambda$ given by \eqref{eq:lambda-principal} specializes to
\begin{equation}
\label{eq:lambda-principal-rank2}
\Lambda =
\left(\!\!\begin{array}{cccc}
0 & 0 & -c & 0 \\
0 & 0 & 0 & -b \\
c & 0 & 0 & bc \\
0 & b & -bc & 0 \\
\end{array}\!\!\right)
\end{equation}

We denote the corresponding quantum cluster algebra by $\AA^\bullet (b,c)$.

The equations \eqref{eq:q-exchange-rel-xx} specialize to
\begin{equation}
\label{eq:q-exchange-rank2}
X'_1 = X^{(-1,c,1,0)} \ + \
 X^{(-1,0,0,0)}, \quad X'_2 = X^{(0,-1,0,1)} \ + \
 X^{(b,-1,0,0)} \, .
\end{equation}

Remembering the definition of $E_a$ (cf. \eqref{eq:standard-monomial-non-normalized}, \eqref{eq:E-a} and
\eqref{eq:nu-of-a}), it is easy to see that in our situation it is given by
\begin{equation}
\label{eq:E-a-principal-rank2}
E_a = v^{- c a_1 a_3 \ - \ b a_2 a_4 \ - \ bc [-a_2]_+ a_3} X^{(0,0,a_3, a_4)} (X'_1)^{[-a_1]_+} X_2^{[a_2]_+} X_1^{[a_1]_+}   (X'_2)^{[-a_2]_+}
\end{equation}
for each $a = (a_1,a_2,a_3,a_4) \in \ZZ^4$.

Turning to the elements $E'_a$, we first note that the mutation $\mu_n = \mu_2$ transforms $\tilde B$ and $\Lambda$ into the matrices
\begin{equation}
\label{eq:tilde-B-prime-principal-rank2}
\tilde B' =
\left(\!\!\begin{array}{cc}
0 & b \\
-c & 0 \\
1 & 0 \\
c & -1 \\
\end{array}\!\!\right), \quad
\Lambda' =
\left(\!\!\begin{array}{cccc}
0 & 0 & -c & 0 \\
0 & 0 & -bc & b \\
c & bc & 0 & bc \\
0 & -b & -bc & 0 \\
\end{array}\!\!\right) \ .
\end{equation}
Then \eqref{eq:Xk-double-prime-expansion} specializes to
\begin{equation}
\label{eq:X1-double-prime-principal-rank2}
 X''_1 = X'_1 (X'_2)^c - \sum_{s=1}^{c} v^{b s^2} \begin{bmatrix} c \\ s \\ \end{bmatrix}_{v^{2b}}
 X^{(bs-1,0,1,c-s)} \, .
 \end{equation}
A direct check shows that, for each $a = (a_1,a_2,a_3,a_4) \in \ZZ^4$, we have
\begin{eqnarray}
\label{eq:E-prime-a-principal-rank2}
& E'_a & = v^{bc ([-a_1]_+ [-a_2]_+  \ + \ [-a_1]_+ a_4 \ - \ c [-a_1]_+ a_3 \ - \ [a_2]_+ a_3 ) \ - \ c a_1 a_3 \ + \ b a_2 a_4}
\\
\nonumber
& & \cdot \
X^{(0,0,a_3, a_4)} X_2^{[-a_2]_+} X_1^{[a_1]_+} (X'_2)^{[a_2]_+} (X''_1)^{[-a_1]_+} \ .
\end{eqnarray}

\begin{proposition}
\label{pr:phi-a-rank2-princial}
Theorem~\ref{th:E-primes-crystal-lattice}  holds for $\AA^\bullet (b,c)$, with the
bijection $\varphi: \ZZ^4 \to \ZZ^4$ given as follows:
\begin{equation}
\label{eq:phi-rank2}
\varphi (a_1,a_2,a_3,a_4) = (a_1, -c [-a_1]_+ - a_2, a_3, a_4 + \min(c[-a_1]_+, [-a_2]_+)) \, .
\end{equation}
\end{proposition}

The main idea of our proof of Proposition~\ref{pr:phi-a-rank2-princial} is to include both bases $\{E_a\}$ and $\{E'_a\}$ into a larger family of ``crystal monomials."
Namely, let $I$ denote the set of $7$-tuples of integers $\mm = (m_3,m_4,m'_1,m_2,m_1, m'_2, m''_1)$ such that the last $5$ components are nonnegative, and let
$I_0 = \{\mm \in I: m'_1 m_1 m''_1 = 0\}$.
For $\mm = (m_3,m_4,m'_1,m_2,m_1, m'_2, m''_1) \in I$, we define
\begin{equation}
\label{eq:M-circ-rank2}
M_\mm^\circ  =  X^{(0,0,m_3,m_4)} (X'_1)^{m'_1} X_2^{m_2} X_1^{m_1} (X'_2)^{m'_2} (X''_1)^{m''_1} \, .
\end{equation}
Then we define the \emph{leading term} $LT(M_\mm^\circ)$ as an element of the ambient quantum torus $\TT$
obtained from $M_{\mm}^\circ$
by replacing the variable $X'_1$ with $X^{(-1,0,0,0)}$, the variable $X'_2$ with $X^{(0,-1,0,1)}$, and the variable $X''_1$ with $X^{(-1,0,1,c)}$.
We set
\begin{equation}
\label{eq:M-rank2}
M_\mm  =  v^{\nu(\mm)} M_{\mm}^\circ \, ,
\end{equation}
where the exponent $\nu(\mm)$ is determined from the following condition:
\begin{equation}
\text{the element $v^{\nu(\mm) -cm'_1 m''_1} LT(M_{\mm}^\circ)$is invariant under the bar-involution.}
\end{equation}
A direct calculation provides the following explicit expression for $\nu(\mm)$:
\begin{align}
\label{eq:nu-m-explicit-rank2}
\nu(\mm) &= c(m'_1 - m_1 - (bc-1)m''_1 - bm'_2) m_3 + b(cm''_1 + m'_2 - m_2) m_4\\
\nonumber
& + c m_1 m''_1 + b m_2 m'_2 + bc m_2 m''_1 \ .
\end{align}

Comparing this definition with the definitions of the elements $E_a$ and $E'_a$, it is easy to check the following: for every
$a = (a_1,a_2, a_3, a_4) \in \ZZ^4$, we have
\begin{equation}
\label{eq:E-M}
E_a = M_{(a_3,a_4, [-a_1]_+, [a_2]_+, [a_1]_+, [-a_2]_+, 0)}, \quad
E'_a = M_{(a_3,a_4, 0, [-a_2]_+, [a_1]_+, [a_2]_+, [-a_1]_+)} \, .
\end{equation}
It follows that the elements $E_a$ (resp.~$E'_a$) are exactly the elements $M_{\mm}$ with
$m'_1 m_1 = m_2 m'_2 = m''_1 = 0$ (resp.~$m'_1 = m_2 m'_2 = m_1 m''_1 = 0$); more precisely, we have
\begin{align}
\label{eq:M-E}
& M_{(m_3,m_4,m'_1,m_2,m_1, m'_2, 0)} = E_{(m_1-m'_1, m_2 - m'_2, m_3,m_4)}\quad (m'_1 m_1 = m_2 m'_2 = 0)\ , \\
\nonumber
& M_{(m_3,m_4,0,m_2,m_1, m'_2,m''_1)} = E'_{(m_1-m''_1, m'_2 - m_2, m_3,m_4)} \quad (m_1 m''_1 = m_2 m'_2 = 0)\ .
\end{align}

Let $\pi: I_0 \to \ZZ^4$ be the mapping given by
\begin{align}
\nonumber
& \pi(m_3,m_4,m'_1,m_2,m_1, m'_2, m''_1) = \\
\label{eq:pi-M-E}
& (m'_1 - m_1 + m''_1, m_2 - m'_2 - c(m''_1 - \underline{m}), \\
\nonumber
& m_3 + \underline{m},
m_4 + \min(m_2 + c \underline{m}, m'_2 + c m''_1)) \ ,
\end{align}
where we abbreviate $\min(m_1, m''_1) = \underline{m}$.
We deduce Proposition~\ref{pr:phi-a-rank2-princial} from the following lemma.

\begin{lemma}
\label{lem:M-crystal}
\begin{enumerate}
\item All elements $M_\mm$ for $\mm \in I$ belong to $\AA_+$.
More precisely, if $\mm \in I_0$ then $M_\mm \in \AA_+ - v \AA_+$, and
if $\mm \in I - I_0$ then $M_\mm \in v \AA_+$.

\item For every $\mm \in I_0$, we have $M_\mm - E_{\pi(\mm)} \in v \AA_+$ .
\end{enumerate}
\end{lemma}

In view of the second equality in \eqref{eq:E-M}, to show that Part (2) of Lemma~\ref{lem:M-crystal} implies Proposition~\ref{pr:phi-a-rank2-princial}
it is enough to prove that
$$\pi(a_3,a_4, 0, [-a_2]_+, [a_1]_+, [a_2]_+, [-a_1]_+) = \varphi(a)$$
for every $a = (a_1,a_2, a_3, a_4) \in \ZZ^4$; but this is immediate from the definitions.

\smallskip

{\bf Proof of Lemma~\ref{lem:M-crystal}.}
First of all, note that \eqref{eq:q-exchange-rank2-product} implies
\begin{equation}
\label{eq:q-exchange-rank2-product}
X'_1 X_1 =  1 \ + \ v^c X^{(0,0,1,0)} X_2^c, \quad X_2 X'_2 = v^{-b} X^{(0,0,0,1)} \ + \
 X_1^b \, .
\end{equation}
In a similar fashion we get
\begin{equation}
\label{eq:q-second-exchange-rank2-product}
X_1 X''_1 =  v^{-c} X^{(0,0,1,c)} \ + \ (X'_2)^c \, .
\end{equation}

In the calculation below we will use the following easily verified property of the ordering $X'_1, X_2, X_1, X'_2, X''_1$
 of cluster variables that appear in the definition \eqref{eq:M-rank2}: any two elements adjacent in this ordering commute with each other.
 Also all these variables quasi-commute with the elements of the form $X^{(0,0,m_2,m_4)}$, with the multiples of quasi-commutation governed by
 the matrix $\Lambda$ in \eqref{eq:lambda-principal-rank2}.
 Using these properties, it is a matter of a routine check to deduce from the equalities \eqref{eq:q-exchange-rank2-product}, \eqref{eq:q-second-exchange-rank2-product} and
 \eqref{eq:X1-double-prime-principal-rank2} the following relations between the elements $M_\mm$.

 Let $\mm = (m_3,m_4,m'_1,m_2,m_1, m'_2, m''_1) \in I$. Then we have:
 \begin{itemize}
 \item if $m'_1 m_1 > 0$ then
 \begin{align}
 \label{eq:M-m'1-m1}
 M_\mm & = v^{cm''_1} M_{(m_3,m_4,m'_1-1,m_2,m_1-1, m'_2, m''_1)}\\
 \nonumber
 & + v^{c(m'_1 + m_1 - 1)} M_{(m_3+1,m_4,m'_1-1,m_2+c,m_1-1, m'_2, m''_1)}\ .
 \end{align}

 \item if $m_2 m'_2 > 0$ then
 \begin{align}
 \label{eq:M-m2-m'2}
 M_\mm & =  M_{(m_3,m_4+1,m'_1,m_2-1,m_1, m'_2-1, m''_1)}\\
 \nonumber
 & + v^{b(m_2 + m'_2 - 1+ cm''_1)} M_{(m_3,m_4,m'_1,m_2-1,m_1+b, m'_2-1, m''_1)}\ .
 \end{align}

 \item if $m_1 m''_1 > 0$ then
 \begin{align}
 \label{eq:M-m1-m''1}
 M_\mm & = v^{cm'_1} M_{(m_3+1,m_4+c,m'_1,m_2,m_1-1, m'_2, m''_1-1)}\\
 \nonumber
 & + v^{c(m_1 + m''_1 - 1)} M_{(m_3,m_4,m'_1,m_2,m_1-1, m'_2+c, m''_1-1)}\ .
 \end{align}

 \item if $m_1 = 0, \, m''_1 > 0$ then
 \begin{align}
 \label{eq:M-m''1}
 M_\mm & = M_{(m_3,m_4,m'_1+1,m_2,0, m'_2+c, m''_1-1)}\\
 \nonumber
 & - \sum_{s=1}^{c} v^{cm'_1 + bs(m_2 + m'_2 + s)} \begin{bmatrix} c \\ s \\ \end{bmatrix}_{v^{2b}}
 M_{(m_3+1,m_4+c-s,m'_1,m_2,bs-1, m'_2, m''_1-1)}\ .
 \end{align}

 \end{itemize}

 Now we are ready to prove Lemma~\ref{lem:M-crystal}.
 We start with Part (1).

 Consider first the special case when $m''_1 = 0$.
 Then $\mm \in I_0$, so we need to prove that $M_\mm \in \AA_+ - v \AA_+$.
 If $m'_1 m_1 = m_2 m'_2 = 0$, then by the first equality in \eqref{eq:M-E}
 $M_\mm$ is one of the elements $E_a$, and there is nothing to prove.
 Thus it suffices to treat the case when either $m'_1 m_1 > 0$, or $m_2 m'_2 > 0$.
 Then $M_\mm$ satisfies (at least) one of the identities  \eqref{eq:M-m'1-m1} and \eqref{eq:M-m2-m'2}.
 In both identities the first term on the right appears with coefficient $1$ (for \eqref{eq:M-m'1-m1} this is our assumption
$m''_1 = 0$), while the coefficient of the second one is a positive power of $v$.
Also both terms on the right have the value of $m'_1 + m'_2$ smaller than that of $\mm$.
Thus, the desired inclusion $M_\mm \in \AA_+ - v \AA_+$ follows by induction on $m'_1 + m'_2$.

It remains to treat the case $m''_1 > 0$.
Then $M_\mm$ satisfies one of the identities \eqref{eq:M-m1-m''1} and \eqref{eq:M-m''1}.
If $\mm \in I - I_0$, that is if $m'_1 m_1 > 0$, then $M_\mm$ satisfies \eqref{eq:M-m1-m''1}, where both coefficients on the right are
positive powers of $v$.
On the other hand, if $\mm \in I_0$, then it can satisfy either of \eqref{eq:M-m1-m''1} and \eqref{eq:M-m''1}, and in each case
the first term on the right appears with coefficient $1$, while all the coefficients in the second one are in $v \ZZ[v]$.
Note also that all terms on the right in each of the identities \eqref{eq:M-m1-m''1} and \eqref{eq:M-m''1} have the value of $m''_1$ smaller than that of $\mm$.
Thus the claims in Part (1) follow by induction on $m''_1$.

To prove Part (2) suppose $\mm \in I_0$.
If $m'_1 m_1 = m_2 m'_2 = m''_1 = 0$ then $M_\mm$ is given by the first equality in \eqref{eq:M-E},
and an easy inspection shows that $M_\mm = E_{\pi(\mm)}$, proving the desired claim.
Thus we can assume that at least one of  $m'_1 m_1, \, m_2 m'_2$ and $m''_1$ is non-zero.
Then $M_\mm$ satisfies at least one of the identities \eqref{eq:M-m'1-m1} - \eqref{eq:M-m''1}.
In all of them the first term on the right appears with coefficient $1$, and we refer to this term as $M_{\mm^-}$.
Furthermore, the rest of the terms on the right appear with coefficients in $v \ZZ[v]$.
Thus Part (1) implies that in each of the cases we have $M_\mm - M_{\mm^-} \in v \AA_+$.
To make the correspondence $\mm \mapsto \mm^-$ well-defined, we take the first term in the first of the identities
\eqref{eq:M-m1-m''1}, \eqref{eq:M-m''1}, \eqref{eq:M-m2-m'2} and \eqref{eq:M-m'1-m1} (in this order!) that is applicable to $M_\mm$.
An easy inspection shows that $\pi(\mm^-) = \pi(\mm)$, and that the iteration of the mapping $\mm \mapsto \mm^-$ terminates (that is, ends up in the set
$\{\mm \in I_0: m'_1 m_1 = m_2 m'_2 = m''_1 = 0\}$) after a finite number of steps for any
initial $\mm \in I_0$.
This completes the proofs of Lemma~\ref{lem:M-crystal} and Proposition~\ref{pr:phi-a-rank2-princial}.


\section{Proof of Theorem~\ref{th:canonical-lusztig-type-basis}-III. Principal quantization in any rank}
\label{sec:proof-main-principal}

In this section we finish the proof of Theorem~\ref{th:canonical-lusztig-type-basis} by dealing with the principal
acyclic quantum cluster algebra of an arbitrary rank~$n$.
Thus, the matrices $\tilde B$ and $\Lambda$ are given by \eqref{tilde-B-principal} and \eqref{eq:lambda-principal},
and we assume that the matrix entries of the exchange matrix~$B$ satisfy \eqref{eq:signs-B-natural-order}.

The arguments below follow those in Section~\ref{sec:proof-main-principal-rank2} with some necessary modifications.
As before, we include both bases $\{E_a\}$ and $\{E'_a\}$ (where now $a$ runs over~$\ZZ^{2n}$) into a larger family of ``crystal monomials."
Namely, let $I$ denote the set of $(4n-1)$-tuples of integers
$$\mm = (m_{n+1}, \dots, m_{2n}, m'_1, \dots, m'_{n-1}, m_1, \dots, m_{n}, m'_n, m''_1, \dots, m''_{n-1})$$
such that the last $3n-1$ components are nonnegative, and
\begin{equation}
\label{eq:m'-m''-restriction}
\text{$m'_i m''_j = 0$ for  $1 \leq j < i \leq n-1$.}
\end{equation}
For $\mm \in I$, we define
\begin{align}
\label{eq:M-circ}
M_\mm^\circ =&  X^{(0,\dots,0,m_{n+1}, \dots, m_{2n})} (X'_1)^{m'_1} \cdots (X'_{n-1})^{m'_{n-1}} \, \times\\
\nonumber
& X^{(m_{1}, \dots, m_{n},0,\dots,0)} (X'_n)^{m'_n} (X''_1)^{m''_1} \cdots (X''_{n-1})^{m''_{n-1}}\, .
\end{align}

Note that the factors in $M_\mm^\circ$ satisfy the following easily checked commutation properties:
\begin{itemize}

\item The elements $X_1, \dots, X_n$ commute with each other.

\item For $i \in [1,n]$ the element $X'_i$ commutes with $X_j$ for $j \in [1,n], j \neq i$.

\item For $i \in [1,n-1]$, the element $X''_i$ commutes with $X'_n$ and with $X_j$ for $j \in [1,n-1], j \neq i$.
\end{itemize}
We will use these properties without further notice.

We define the \emph{leading term} $LT(M_\mm^\circ)$ as an element of the ambient quantum torus $\TT$
obtained from $M_{\mm}^\circ$ by the following replacement (for $j \in [1,n-1]$):
\begin{equation}
\label{eq:LT-principal}
X'_j \mapsto X^{- e_j - b_j^{<j}}, \quad X'_n \mapsto X^{-e_n + e_{2n}}, \quad
X''_j \mapsto X^{- e_j + b_j^{>j} +b_{nj}(e_{2n}-e_n)} \ ,
\end{equation}
where as before, $b_j$ stands for the $j$th column of $\tilde B$, and $b_j^{<j}$ and $b_j^{>j}$ are the corresponding
truncated vectors.
We abbreviate
\begin{equation}
\label{eq:nu'}
\nu'(\mm) = \sum_{j=1}^{n-1} d_j m'_j m''_j \ ,
\end{equation}
and we set
\begin{equation}
\label{eq:M}
M_\mm  =  v^{\nu(\mm)} M_{\mm}^\circ \, ,
\end{equation}
where the exponent $\nu(\mm)$ is determined from the following condition:
\begin{equation}
\label{eq:M-normalization-anyrank}
\text{the element $v^{\nu(\mm)-\nu'(\mm)} LT(M_{\mm}^\circ)$ is invariant under the bar-involution.}
\end{equation}

As in the case $n=2$, it is easy to check that both triangular bases $\{E_a\}$ and $\{E'_a\}$ are contained in the family of elements $M_\mm$.
More precisely, the elements $E_a$  are exactly the elements $M_{\mm}$ with
\begin{equation}
\label{eq:M-E-conditions}
m'_1 m_1 = m'_2 m_2 = \cdots = m'_n m_n = m''_1 = m''_2 = \cdots = m''_{n-1} =  0 \ .
\end{equation}
Namely, for every $\mm \in I$ satisfying \eqref{eq:M-E-conditions}, we have $M_\mm = E_a$, where the vector $a \in \ZZ^{2n}$ is related to $\mm$ as follows:
\begin{align}
\label{eq:M-E-any-rank}
& a_i = m_i \quad (i \in [n+1,2n]), \,\, a_j = m_j - m'_j \quad (j \in [1,n]) \ ,\\
\nonumber
& m_i = a_i \quad (i \in [n+1,2n]), \,\, m_j = [a_j]_+, \,\, m'_j = [-a_j]_+ \quad (j \in [1,n]) \ .
\end{align}
Similarly, the elements $E'_a$  are exactly the elements $M_{\mm}$ with
\begin{equation}
\label{eq:M-E'-conditions}
m_1 m''_1 = m_2 m''_2 = \cdots = m_{n-1} m''_{n-1} = m'_n m_n = m'_1 = m'_2 = \cdots = m'_{n-1} =  0 \ ,
\end{equation}
and for every $\mm \in I$ satisfying \eqref{eq:M-E'-conditions}, we have $M_\mm = E'_a$, where the vector $a \in \ZZ^{2n}$ is related to $\mm$ as follows:
\begin{align}
\nonumber
& a_i = m_i \quad (i \in [n+1,2n]), \,\, a_j = m_j - m'_j \quad (j \in [1,n-1]), \,\, a_n = m'_n - m_n  \ ,\\
\label{eq:M-E'-any-rank}
& m_i = a_i \quad (i \in [n+1,2n]), \,\, m_j = [a_j]_+, \,\, m'_j = [-a_j]_+ \quad (j \in [1,n-1]) \ ,\\
\nonumber
& m_n = [-a_n]_+, \,\, m'_n = [a_n]_+\ .
\end{align}

As in the previous section, to finish the proof of Theorem~\ref{th:canonical-lusztig-type-basis} it suffices to prove the following
generalization of Lemma~\ref{lem:M-crystal}.

\begin{proposition}
\label{pr:M-crystal-anyrank}
Let
$I_0 = \{\mm \in I: m'_j, m_j, m''_j = 0 \,\, {\rm for}\,\, j \in [1,n-1]\}$.
\begin{enumerate}
\item All elements $M_\mm$ for $\mm \in I$ belong to $\AA_+$.
More precisely, if $\mm \in I_0$ then $M_\mm \in \AA_+ - v \AA_+$, and
if $\mm \in I - I_0$ then $M_\mm \in v \AA_+$.

\item There exists a mapping $\pi: I_0 \to \ZZ^{2n}$ such that, for every $\mm \in I_0$, we have $M_\mm - E_{\pi(\mm)} \in v \AA_+$.
\end{enumerate}
\end{proposition}

The proof basically follows that of Lemma~\ref{lem:M-crystal} but with some modifications.
First we introduce a subset $I_{00} \subset I_0$ by setting
$$I_{00} =  \{\mm \in I: m''_j = 0 \,\, {\rm for}\,\, j \in [1,n-1]\} \ .$$
Let $\tilde \AA_+$ denote the $\ZZ[v]$-submodule of $\AA$ generated by the $M_\mm$ for $\mm \in I_{00}$.
Clearly, Proposition~\ref{pr:M-crystal-anyrank} is a consequence of the following two lemmas.

\begin{lemma}
\label{lem:killing-m''}
\begin{enumerate}
\item All elements $M_\mm$ for $\mm \in I - I_{00}$ belong to $\tilde \AA_+$.
\item If $\mm \in I - I_0$ then $M_\mm \in v \tilde \AA_+$.
\item There exists a mapping $\pi_0: I_0 \to I_{00}$ such that, for every $\mm \in I_0$, we have $M_\mm - M_{\pi_0(\mm)} \in v \tilde \AA_+$.
\end{enumerate}
\end{lemma}

\begin{lemma}
\label{lem:killing-m-m'}
 There exists a mapping $\pi_{00}: I_{00} \to \ZZ^{2n}$ such that, for every $\mm \in I_{00}$, we have $M_\mm - E_{\pi_{00}(\mm)} \in v \AA_+$.
In particular, we have $\tilde \AA_+ = \AA_+$.
\end{lemma}

Proceeding by induction on $\sum_{j=1}^{n-1} m''_j$, we see that Lemma~\ref{lem:killing-m''} is a consequence of the following two identities.

\begin{lemma}
\label{lem:identities-m''}
Suppose $\mm \in I$ is such that $\sum_{j=1}^{n-1} m''_j > 0$, and let $j \in [1,n-1]$ be the smallest index such that  $m''_j > 0$.
\begin{enumerate}
\item If $m_j > 0$ then we have
\begin{equation}
\label{eq:identity-m-m''}
M_\mm = v^{d_j m'_j} M_{\mm^+} + v^{d_j (m_j + m''_j - 1)} M_{\mm^-} \ ,
\end{equation}
where $\mm^+$ is obtained from $\mm$ by the replacement
\begin{align*}
& m_j \mapsto m_j -1, \, m''_j \mapsto m''_j - 1, \, m_k \mapsto m_k + b_{kj} \,\, (j < k \leq n-1) \ , \\
& m_{n+j} \mapsto m_{n+j} + 1, \, m_{2n} \mapsto m_{2n} + b_{nj} \, ,
\end{align*}
while $\mm^-$ is obtained from $\mm$ by the replacement
$$m_j \mapsto m_j -1, \, m''_j \mapsto m''_j - 1, \, m_i \mapsto m_i - b_{ij} \,\, (1 \leq i < j), \,
m'_{n} \mapsto m'_{n} + b_{nj} \, .$$
\item If $m_j=0$ then we have
\begin{equation}
 \label{eq:identity-m''}
 M_\mm  = M_{\mm^+} - \sum_{s=1}^{b_{nj}} v^{d_j m'_j + d_n s (m_n + m'_n + s)} \begin{bmatrix} b_{nj} \\ s \\ \end{bmatrix}_{v^{2d_n}}
 M_{\mm^-(s)}\ ,
 \end{equation}
 where $\mm^+$ is obtained from $\mm$ by the replacement
 $$m''_j \mapsto m''_j - 1, \, m'_j \mapsto m'_j + 1, \,
m'_{n} \mapsto m'_{n} + b_{nj} \, ,$$
and $\mm^-(s)$ is obtained from $\mm$ by the replacement
 \begin{align*}
& m''_j \mapsto m''_j - 1, \, m_j \mapsto m_j - sb_{jn} - 1, \,
m_i \mapsto m_i - s b_{in} \,\, (1 \leq i < j), \\
& m_k \mapsto m_k + b_{kj} - s b_{kn} \,\, (j < k \leq n-1),
m_{n+j} \mapsto m_{n+j} + 1, \\
& m_{2n} \mapsto m_{2n} + b_{nj} - s \, .
\end{align*}
\end{enumerate}
\end{lemma}

\noindent {\bf Proof of Lemma~\ref{lem:identities-m''}.} To prove \eqref{eq:identity-m-m''}, we use the identity
\begin{equation}
\label{eq:Xj-Xj''}
X_j X''_j = v^{-d_j} X^{b_j^{>j} + b_{nj} (e_{2n}-e_n)} + X^{-b_j^{<j}} (X'_n)^{b_{nj}} \quad (1 \leq j < n) \
\end{equation}
which is an easy consequence of \eqref{eq:X-double-prime-k}.
Recalling the definition of $M_\mm^\circ$ in \eqref{eq:M-circ}, and the commutation properties collected after it, we can express the factor
$X^{(m_{1}, \dots, m_{n},0,\dots,0)}$ in $M_\mm^\circ$  as $X^{(m_{1}, \dots, m_{n},0,\dots,0) - m_j e_j} X_j^{m_j}$, and then interchange the term
$(X''_j)^{m''_j}$ with the commuting term $(X'_n)^{m'_n}$, so that it will stand right after $X_j^{m_j}$.
Using \eqref{eq:Xj-Xj''}, we replace  $X_j^{m_j} (X''_j)^{m''_j}$ with the sum of two terms: $v^{-d_j} X_j^{m_j-1} X^{b_j^{>j} + b_{nj} e_{2n}}(X''_j)^{m''_j-1}$, and
$X_j^{m_j-1}X^{-b_j^{<j}} (X'_n)^{b_{nj}} (X''_j)^{m''_j-1}$.
Again using the commutation relations we conclude that
$$M_\mm^\circ = v^{-d_j} M_{\mm^+}^\circ + M_{\mm^-}^\circ \ .$$
In view of \eqref{eq:M}, this implies the following:
\begin{equation}
\label{eq:identity-m-m''-thru-nu}
M_\mm = v^{\nu(\mm)-\nu(\mm^+)-d_j} M_{\mm^+} + v^{\nu(\mm)-\nu(\mm^-)} M_{\mm^-} \ .
\end{equation}
Here the exponents of $v$ are determined from the conditions that each of the elements
$v^{\nu(\mm)-\nu'(\mm)} LT(M_{\mm}^\circ)$, $v^{\nu(\mm^+)-\nu'(\mm^+)} LT(M_{\mm^+}^\circ)$,
and $v^{\nu(\mm^-)-\nu'(\mm^-)} LT(M_{\mm^-}^\circ)$ is invariant under the bar-involution (see \eqref{eq:M-normalization-anyrank});
recall that the function $\nu'(\mm)$ is defined by \eqref{eq:nu'}.

Remembering the definition \eqref{eq:LT-principal} of the leading term, we see that $LT(M_{\mm^+}^\circ)$
is obtained from $LT(M_\mm^\circ)$ by replacing $LT(X_j X''_j) = X_j X^{-e_j + b_j^{>j} + b_{nj} (e_{2n}-e_n)}$ with
$X^{b_j^{>j} + b_{nj} (e_{2n}-e_n)}$ in the appropriate place in the product expansion of $LT(M_\mm^\circ)$.
Since
$$X_j X^{-e_j + b_j^{>j} + b_{nj} (e_{2n}-e_n)} = v^{-d_j} X^{b_j^{>j} + b_{nj} (e_{2n}-e_n)}$$
(see \eqref{eq:multiplication-quantum-torus} and \eqref{eq:lambda-principal}), we see that
$LT(M_{\mm^+}^\circ) = v^{d_j} LT(M_{\mm}^\circ)$.
This implies that $\nu(\mm)-\nu'(\mm) = d_j + \nu(\mm^+)-\nu'(\mm^+)$.
Therefore, we have
$$\nu(\mm)-\nu(\mm^+)-d_j = \nu'(\mm) -  \nu'(\mm^+) =  d_j m'_j m''_j - d_j m'_j (m''_j - 1) = d_j m'_j$$
showing that $M_{\mm^+}$ appears in the right hand side of \eqref{eq:identity-m-m''-thru-nu} with the same coefficient
$v^{d_j m'_j}$ as in the right hand side of \eqref{eq:identity-m-m''}.

To prove the similar statement for  $M_{\mm^-}$, we note that the elements $LT(M_{\mm}^\circ)$ and $LT(M_{\mm^-}^\circ)$
can be factored as follows:
$$LT(M_{\mm}^\circ) = v^{\gamma - d_j} X^f X^{b_j^{>j} + b_{nj} (e_{2n}-e_n)} X^g, \quad
LT(M_{\mm^-}^\circ) = v^{\gamma} X^f X^{-b_j^{<j} + b_{nj} (e_{2n}-e_n)} X^g \ ,$$
where $\gamma \in \ZZ$, and $f, g \in \ZZ^{2n}$ are integer vectors such that $e_j$ appears in $f$ with the coefficient
$m_j - m'_j - 1$ and appears in $g$ with the coefficient $-(m''_j - 1)$ (the statement about the coefficient of $e_j$ in $f$ uses the fact
that, according to \eqref{eq:m'-m''-restriction}, we have $m'_i = 0$ for  $1 \leq j < i \leq n-1$).
Using \eqref{eq:multiplication-quantum-torus} and the compatibility condition \eqref{eq:orthogonality}, we conclude that
$\nu(\mm)-\nu'(\mm) = d_j + \nu(\mm^-)-\nu'(\mm^-) + d_j (m_j - m'_j + m''_j - 2)$.
It follows that
$$\nu(\mm)-\nu(\mm^-) = \nu'(\mm) - \nu'(\mm^-) + d_j (m_j - m'_j + m''_j - 1) =  d_j (m_j + m''_j - 1) \ .$$
Thus the desired identity \eqref{eq:identity-m-m''} becomes a consequence of \eqref{eq:identity-m-m''-thru-nu}, finishing the proof of Part~1 of
Lemma~\ref{lem:identities-m''}.

The identity \eqref{eq:identity-m''} is proved in a similar way, and we leave the details to the reader (our starting point is the identity
\begin{equation}
\label{eq:X''j-principal}
 X''_j = X'_j (X'_n)^{b_{nj}} - \sum_{s=1}^{b_{nj}} v^{s^2 d_n} \begin{bmatrix} b_{nj} \\ s \\ \end{bmatrix}_{v^{2d_n}}
 X^{-e_j + b_j^{>j} + b_{nj}(e_{2n} - e_n) - sb_n}
 \end{equation}
which is a specialization of \eqref{eq:Xk-double-prime-expansion}).
This completes the proof of Lemma~\ref{lem:identities-m''}.

\smallskip

Now we turn to the proof of Lemma~\ref{lem:killing-m-m'}.
We view the index set $I_{00}$ as the set of
$3n$-tuples of integers
$$\mm = (m_{n+1}, \dots, m_{2n}, m'_1, \dots, m'_{n}, m_1, \dots, m_{n})$$
such that the last $2n$ components are nonnegative.
Once again we include the family of monomials $\{M_\mm: \mm \in I_{00}\}$ into a larger family.
Namely, we define
$$\hat I_{00} = I_{00} \times [1,n-1] \ ,$$
and for each $(\mm, j) \in \hat I_{00}$, define an element $M^\circ_{\mm;j}$ by setting
\begin{align}
\label{eq:Mj-circ}
M_{\mm;j}^\circ =&  X^{(0,\dots,0,m_{n+1}, \dots, m_{2n})} (X'_1)^{m'_1} \cdots (X'_{j})^{m'_{j}} \, \times\\
\nonumber
& X^{(m_{1}, \dots, m_{n},0,\dots,0)} (X'_{j+1})^{m'_{j+1}} \cdots (X'_{n})^{m'_{n}}\, .
\end{align}
For $(\mm, j) \in \hat I_{00}$, we define the \emph{leading term} $LT_j (M^\circ_{\mm;j})$ as an element of the quantum torus $\TT$
obtained from $M_{\mm;j}^\circ$ by the following replacement:
\begin{equation}
\label{eq:LT-j}
X'_i \mapsto X^{- e_i - b_i^{<i}}   \quad (1 \leq i \leq j), \quad X'_k \mapsto X^{- e_k + b_k^{>k}}   \quad (j+1 \leq k \leq n) \ .
\end{equation}
Finally we set
\begin{equation}
\label{eq:Mj}
M_{\mm;j}  =  v^{\nu(\mm;j)} M_{\mm;j}^\circ \, ,
\end{equation}
where the exponent $\nu(\mm;j)$ is determined from the following condition:
\begin{equation}
\label{eq:M-j-normalization}
\text{the element $v^{\nu(\mm;j)} LT_j(M_{\mm;j}^\circ)$ is invariant under the bar-involution.}
\end{equation}

Comparing this definition with \eqref{eq:M-circ} and \eqref{eq:M-normalization-anyrank}, we see that $M_{\mm}^\circ = M_{\mm;n-1}^\circ$
and $M_{\mm} = M_{\mm;n-1}$ for every $\mm \in I_{00}$.
Thus, to prove Lemma~\ref{lem:killing-m-m'}, it suffices to show the following:
\begin{equation}
\label{eq:pi-hat}
\text{For every $(\mm,j) \in \hat I_{00}$ there exists $a \in \ZZ^{2n}$ such that
$M_{\mm;j} - E_a \in v \AA_+$}.
\end{equation}

We deduce \eqref{eq:pi-hat} from the following identities.

\begin{lemma}
\label{lem:identities-m'-m}
Suppose $\mm \in I_{00}$ is such that $m'_j m_j > 0$ for some $j \in [1,n]$.
For $\mm \in I_{00}$, let $\mm[<j] \in I_{00}$ denote an element obtained from $\mm$ by the replacement
\begin{equation*}
m_j \mapsto m_j -1, \, m'_j \mapsto m'_j - 1, \, m_i \mapsto m_i - b_{ij} \,\, (1 \leq i < j) \ ,
\end{equation*}
and $\mm[>j] \in I_{00}$ denote an element obtained from $\mm$ by the replacement
$$m_j \mapsto m_j -1, \, m'_j \mapsto m'_j - 1, \, m_k \mapsto m_k + b_{kj} \,\, (j < k \leq n), \,
m_{n+j} \mapsto m_{n+j} + 1 \, .$$
\begin{enumerate}
\item If $j < n$ then we have
\begin{equation}
\label{eq:identity-m'-m}
M_{\mm;j} = M_{\mm[<j];j} + v^{d_j (m_j + m'_j - 1)} M_{\mm[>j];j} \ .
\end{equation}
\item If $1 < j$ then we have
\begin{equation}
\label{eq:identity-m-m'}
M_{\mm;j-1} = M_{\mm[>j];j-1} + v^{d_j (m_j + m'_j - 1)} M_{\mm[<j];j-1} \ .
\end{equation}
\end{enumerate}
\end{lemma}

The identities \eqref{eq:identity-m'-m} and \eqref{eq:identity-m-m'} are proved in the same way as \eqref{eq:identity-m-m''} above,
with the role of \eqref{eq:Xj-Xj''} played by the identities
\begin{equation}
\label{eq:X'j-Xj}
X'_j X_j =  X^{-b_j^{<j}} + v^{d_j} X^{b_j^{>j}} \quad (1 \leq j < n)
\end{equation}
and
\begin{equation}
\label{eq:Xj-X'j}
X_j X'_j =  v^{-d_j} X^{b_j^{>j}} +  X^{- b_j^{<j}} \quad (1 < j \leq n) \ .
\end{equation}
We leave the details to the reader.

To deduce \eqref{eq:pi-hat} from Lemma~\ref{lem:identities-m'-m}, we first make an easy observation:
\begin{equation}
\label{eq:Mj-j-1}
\text{If $(\mm,j) \in \hat I_{00}$ is such that $1 < j$, and $m'_j m_j = 0$ then $M_{\mm;j} = M_{\mm;j-1}$.}
\end{equation}
Now let $(\mm,k) \in \hat I_{00}$ be such that $m'_i m_i > 0$ for some $i \in [1,n]$.
Using \eqref{eq:Mj-j-1} if necessary, we can find $j \in [1,n]$ such that $m'_j m_j > 0$, and
$M_{\mm;k}$ is equal either to $M_{\mm;j}$ or to $M_{\mm;j-1}$, so that it satisfies one of the identities
\eqref{eq:identity-m'-m} and \eqref{eq:identity-m-m'}.
The fact that $M_{\mm;k}$  satisfies the desired property \eqref{eq:pi-hat}, then follows by induction on
$r(\mm) = \sum_{i=1}^n m'_i$ (since we have $r(\mm[<j]) = r(\mm[>j]) = r(\mm) - 1$).

This concludes the proof of \eqref{eq:pi-hat} and hence that of Lemma~\ref{lem:killing-m-m'}, Proposition~\ref{pr:M-crystal-anyrank},
and Theorem~\ref{th:canonical-lusztig-type-basis}.

\section{Example: coefficient-free type $A_1^{(1)}$}
\label{sec:coef-free-rank2-affine}

In this section we illustrate the above results by computing the canonical triangular basis in the quantum cluster algebra
$\AA$ with the initial quantum seed given as follows:
\begin{itemize}
\item $m = n = 2, \, d_1 = d_2 = 2$.

\item $\tilde B = B =
\left(\!\!\begin{array}{cc}
0 & -2 \\
2 & 0  \\
\end{array}\!\!\right) \ ,
\Lambda = \left(\!\!\begin{array}{cc}
0 & -1 \\
1 & 0  \\
\end{array}\!\!\right) \ .$

\item $\tilde \XX = \XX = (X_1, X_2)$.
\end{itemize}

This cluster algebra and various bases in it were studied in detail in \cite{ding-xu,lampe}.
In what follows  we use the results in \cite{ding-xu}.
Since our choice of $\tilde B$ and $\Lambda$ differs from that of \cite{ding-xu} by the sign, to reconcile our setups their
formal variable $q$ and our formal variable $v$ are related by $q = v^{-2}$.

As in  \cite{ding-xu}, the set of cluster variables in $\AA$ is denoted by $\{X_m: m \in \ZZ\}$, with the clusters
$\{X_m, X_{m+1}\}$ for $m \in \ZZ$.
The commutation relations in $\AA$ are
\begin{equation}
\label{eq:commutation-A11}
X_{m+1} X_m = v^2 X_m X_{m+1} \quad (m \in \ZZ) \ ,
\end{equation}
the exchange relations are
\begin{equation}
\label{eq:exchange-A11}
X_{m+1} X_{m-1} = v^2 X_m^2 + 1 \quad (m \in \ZZ) \ ,
\end{equation}
and the (bar-involution invariant) cluster monomials are the elements
\begin{equation}
\label{eq:cluster-mon-A11}
v^{a_1 a_2} X_m^{a_1} X_{m+1}^{a_2} \quad (a_1, a_2 \in \ZZ_{\geq 0}, \,\, m \in \ZZ) \ .
\end{equation}

Following \cite{ding-xu} we denote by $X_\delta$ the element of $\AA$ given by
\begin{equation}
\label{eq:Xdelta-A11}
X_\delta = v X_3 X_0 - v^3 X_2 X_1 \ .
\end{equation}
Let $S_{-1}(z), S_0 (z), S_1(z), \dots$ be the sequence of
(normalized) \emph{Chebyshev polynomials of second kind} given by
the initial conditions
\begin{equation}
\label{eq:Cheb-2-initial}
S_{-1}(z) = 0, \quad S_0(z) = 1 \, ,
\end{equation}
and the recurrence relation
\begin{equation}
\label{eq:Cheb-2-recurrence}
S_r(z) = z S_{r-1}(z) -  S_{r-2}(z) \quad (r \geq 1) \, .
\end{equation}
Thus we have $S_1(z) = z, \ S_2(z) = z^2-1$, etc.

\begin{proposition}
\label{pr:can-triangular-basis-A11}
The canonical triangular basis in $\AA$ consists of all cluster monomials given by \eqref{eq:cluster-mon-A11}, and the elements
$S_r(X_\delta)$ for all $r \geq 1$.
\end{proposition}

\begin{remark}
As shown in \cite{ding-xu}, the cluster variables and the elements $S_r(X_\delta)$ for all $r \geq 1$ are precisely the quantum Caldero-Chapoton characters
associated with indecomposable representations of the Kronecker quiver (see \cite{rupel}). In this case the canonial triangular basis coincides with the natural quantum version of the dual semicanonical basis introduced for the commutative setting in \cite{CZ}; this version was also discovered and studied in \cite{lampe}.
\end{remark}

\noindent {\bf Proof of Proposition~\ref{pr:can-triangular-basis-A11}.}
The fact that all cluster monomials belong to the canonical triangular basis $\BB$ in $\AA$, follows from
Corollary~\ref{cor:B-contains-acyclic-cluster-monomials} (in our situation, \emph{all} the seeds are acyclic).
They are labeled by the lattice $\ZZ^2$, and our first task is to describe this labeling explicitly.

An easy calculation shows that the elements of the ``standard" basis $\{E_a: a \in \ZZ^2\}$ can be expressed as follows:
\begin{equation}
\label{eq:E-a-A11}
E_a = v^{a_1 a_2} X_3^{[-a_1]_+} X_1^{[a_1]_+} X_2^{[a_2]_+} X_0^{[-a_2]_+} \ .
\end{equation}
Because of the translational symmetry of the relations \eqref{eq:commutation-A11} and \eqref{eq:exchange-A11}, there are well-defined
mutually inverse algebra automorphisms $\eta_+$ and $\eta_-$ of $\AA$ acting on the generators (that is, the cluster variables) by
\begin{equation}
\label{eq:tau-pm-A11}
\eta_+ (X_m) = X_{m+1}, \,\,  \eta_- (X_m) = X_{m-1} \quad (m \in \ZZ) \ .
\end{equation}
Comparing \eqref{eq:E-a-A11} with \eqref{eq:E-prime-a-principal-rank2} (and ignoring frozen variables), we see
that $\eta_-(E_{(a_1,a_2)}) = E'_{(a_2, a_1)}$.
Thus,  Proposition~\ref{pr:phi-a-rank2-princial} implies that
\begin{equation}
\label{eq:etaEa-Ea-A11}
\eta_-(E_{(a_1,a_2)}) - E_{(a_2, -2 [-a_2]_+ - a_1)} \in v \AA_+ \ .
\end{equation}
Since $\eta_+ = \eta_-^{-1}$, we also have
\begin{equation}
\label{eq:eta+Ea-Ea-A11}
\eta_+(E_{(a_1,a_2)}) - E_{(-2 [-a_1]_+ - a_2, a_1)} \in v \AA_+ \ .
\end{equation}
As a consequence, the elements $C_a$ of $\BB$ satisfy
\begin{equation}
\label{eq:etaCa-A11}
\eta_-(C_{(a_1,a_2)}) = C_{(a_2, -2 [-a_2]_+ - a_1)}, \quad
\eta_+(C_{(a_1,a_2)}) = C_{(-2 [-a_1]_+ - a_2, a_1)} \ .
\end{equation}
Iterating \eqref{eq:etaCa-A11}, we conclude that the cluster monomials in \eqref{eq:cluster-mon-A11} are labeled as follows:
\begin{equation}
\label{eq:cluster-mon-label-A11}
v^{a_1 a_2} X_m^{a_1} X_{m+1}^{a_2} = C_{a_1 \alpha(m) + a_2 \alpha(m+1)} \quad (a_1, a_2 \in \ZZ_{\geq 0}, \,\, m \in \ZZ) \ ,
\end{equation}
where the vectors $\alpha(m) \in \ZZ^2$ are given by
\begin{equation}
\label{eq:cluster-variables-labels-A11}
\alpha(1-r) = (1-r,-r), \quad \alpha(2+r) = (-r,1-r) \quad (r \geq 0) \ .
\end{equation}

As an easy consequence of \eqref{eq:cluster-mon-label-A11} and \eqref{eq:cluster-variables-labels-A11}, we conclude that the cluster monomials are all the elements
$C_a$ for $a \in \ZZ^2 - \{(-r,-r) : r \geq 1\}$.
To finish the proof of Proposition~\ref{pr:can-triangular-basis-A11}, it remains to show that
\begin{equation}
\label{eq:C-rdelta-A11}
 C_{(-r,-r)} = S_r(X_\delta) \quad (r \geq 1) \ .
\end{equation}
We use the following properties of the elements $S_r(X_\delta)$ (all of them are established in \cite{ding-xu}).

\begin{lemma}
\label{lem:SXdelta-A11}
\begin{enumerate}
\item The elements $S_r(X_\delta)$ are invariant under the bar-involution, and under each of the automorphisms
$\eta_+$ and $\eta_-$.
\item For each $r \geq -1$ we have
\begin{equation}
\label{eq:SrXdelta-A11}
S_r(X_\delta) = v^r X_{r+2} X_0 - v^{r+2} X_{r+1} X_1 \ .
\end{equation}
\end{enumerate}
\end{lemma}

Recall the notation $\AA_+$ for the $\ZZ[v]$-linear span of the basis $\{E_a : a \in \ZZ^2\}$.
In this notation, the condition \eqref{eq:b-triangularity} takes the form $C_a - E_a \in v \AA_+$.
Recall that we also have a stronger condition \eqref{eq:b-double-triangularity}.
We have defined the partial order $\prec$ in \eqref{eq:partial-order-on-Zm}.
However, in our current situation it is possible to replace this partial order by the following sharper one:
\begin{equation}
\label{eq:partial-order-on-Z2-A11}
a' = (a'_1, a'_2) \prec a = (a_1,a_2) \Longleftrightarrow [-a'_1]_+ < [-a_1]_+, \,\, [-a'_2]_+ < [-a_2]_+  \, ;
\end{equation}
indeed a direct check shows that the condition \eqref{eq:bar-triangularity} holds for this sharper partial order.

The last ingredient we need to prove \eqref{eq:C-rdelta-A11} is the following lemma.

\begin{lemma}
\label{lem:Ea-times X0-A11}
For every $a = (a_1, a_2) \in \ZZ^2$, we have
\begin{equation}
\label{eq:Ea-X0-A11}
v^{-a_1} E_a X_0 - E_{(a_1, a_2-1)} \in v \AA_+ \ .
\end{equation}
\end{lemma}

\begin{proof} We prove \eqref{eq:Ea-X0-A11} by a direct calculation of the left hand side using
\eqref{eq:E-a-A11}, \eqref{eq:commutation-A11} and \eqref{eq:exchange-A11}.
There are four cases to consider.
In each case we just give the result of a calculation leaving the details to the reader.

\smallskip

\noindent {\bf Case 1:} $a_2 \leq 0$. Then we have $v^{-a_1} E_a X_0 - E_{(a_1, a_2-1)} = 0$.

\smallskip

\noindent {\bf Case 2:} $a_2 > 0, \, a_1 \geq 0$. Then
$v^{-a_1} E_a X_0 - E_{(a_1, a_2-1)} = v^{2a_2} E_{(a_1+2, a_2 - 1)}$.

\smallskip

\noindent {\bf Case 3:} $a_2 > 0, \, a_1 = -1$. Then
$$v^{-a_1} E_a X_0 - E_{(a_1, a_2-1)} = v^{2a_2} E_{(1, a_2 - 1)} + v^{2(a_2 + 2)} E_{(1, a_2 + 1)} \, .$$

\smallskip

\noindent {\bf Case 4:} $a_2 > 0, \, a_1 \leq -2$. Then
\begin{align*}
v^{-a_1} E_a X_0 - E_{(a_1, a_2-1)} &= v^{2a_2} E_{(a_1 + 2, a_2 - 1)} +
(v^{2(a_2 - a_1 - 1)} + v^{2(a_2 - a_1 + 1)}) E_{(a_1+2, a_2 + 1)}\\
&+ v^{2(a_2 - 2a_1)} E_{(a_1+2, a_2 + 3)} \, .
\end{align*}
Since in all the cases the right hand side belongs to $v \AA_+$, we are done.
\end{proof}

Now everything is ready for the proof of \eqref{eq:C-rdelta-A11}.
For $r = 1$, we have
$$S_r(X_\delta) = X_\delta = v X_3 X_0 - v^3 X_2 X_1 = E_{(-1,-1)} - v^4 E_{(1,1)}
\in E_{(-1,-1)} + v \AA_+ \ ;$$
since $X_\delta$ is also invariant under the bar-involution (see Lemma~\ref{lem:SXdelta-A11}), it follows
that $X_\delta = C_{-1,-1}$.
Thus we assume that $r \geq 2$.

As a special case of \eqref{eq:cluster-mon-label-A11}, we have $X_{r+2} = C_{(-r,1-r)}$.
Applying \eqref{eq:b-double-triangularity} with the partial order $\prec$ given by \eqref{eq:partial-order-on-Z2-A11}, we
have
$$X_{r+2} = E_{(-r,1-r)} + \sum_{[-a_1]_+ < r, \ [-a_2]_+ < r - 1} c_{a_1, a_2} E_{(a_1,a_2)} \, $$
where $c_{a_1, a_2} \in v ZZ[v]$ for all $(a_1,a_2)$.
Multiplying both sides on the left with $v^r$ and on the right with $X_0$, and using \eqref{eq:Ea-X0-A11}, we get
$$v^r X_{r+2} X_0 \in E_{(-r,-r)} + \sum_{[-a_1]_+ < r, \ [-a_2]_+ < r - 1}  v^{r + a_1} c_{a_1, a_2} E_{(a_1,a_2-1)}
+ v \AA_+ \subseteq E_{(-r,-r)} + v \AA_+$$
(the last inclusion follows since $r + a_1 > 0$ for $[-a_1]_+ < r$).

By the same token, we have $v^{r-2} X_{r} X_0 \in E_{(2-r,2-r)} + v \AA_+$, implying that
$$v^{r+2} X_{r+1} X_1 = v^4 \eta_+ (v^{r-2} X_{r} X_0) \in v^4  E_{(2-r,2-r)} + v \AA_+ =
v \AA_+ \, .$$
In view of \eqref{eq:SrXdelta-A11}, we conclude that
$$S_r(X_\delta) = v^r X_{r+2} X_0 - v^{r+2} X_{r+1} X_1 \in E_{(-r,-r)} + v \AA_+ \ .$$
Since $S_r(X_\delta)$ is also invariant under the bar-involution (see Lemma~\ref{lem:SXdelta-A11}), it follows
that $S_r(X_\delta) = C_{-r,-r}$, finishing the proofs of \eqref{eq:C-rdelta-A11} and Proposition~\ref{pr:can-triangular-basis-A11}.

\begin{remark}
A big part of the above proof carries over without difficulty to the general rank~$2$ case,  where 
$\tilde B = B =
\left(\!\!\begin{array}{cc}
0 & -b \\
c & 0  \\
\end{array}\!\!\right)$, with arbitrary positive integers $b$ and $c$. 
However the problem of finding explicit expressions for all the elements $C_a$ is still open in this generality.
It definitely deserves a further study.
\end{remark}

\section{Proof of Lusztig's Lemma}
\label{sec:proof-Lusztig-lemma}

In this section we prove Theorem~\ref{th:Lusztig-lemma}.

\begin{proof}
According to \eqref{eq:bar-triangularity}, we have
\begin{equation}
\label{eq:matrix-bar}
\overline {E_a}= E_a + \sum_{a'\in L} r_{a,a'} E_{a'},
\end{equation}
where $r_{a,a'}\in \ZZ[v,v^{-1}]$, and $r_{a,a'}=0$ unless $a' \prec a$.
Expanding the equality $\overline {\overline {E_a}} = E_a$, we see
that the condition that $x \mapsto \overline x$ is an involution
is equivalent to the following: for all $a, a' \in L$, we have
\begin{equation}
\label{eq:r-condition}
r_{a,a'} + \overline{r_{a,a'}} + \sum_{a''} \overline{r_{a,a''}} r_{a'',a'} = 0.
\end{equation}
Applying the bar-involution on both sides
of \eqref{eq:r-condition}, we also get
\begin{equation}
\label{eq:r-condition-twisted}
r_{a,a'} + \overline{r_{a,a'}} + \sum_{a''} r_{a,a''}\overline{r_{a'',a'}} = 0.
\end{equation}

According to \eqref{eq:b-triangularity}, the desired element $C_a$
must have the form
\begin{equation}
\label{eq:matrix-b-e}
C_a= E_a + \sum_{a'\in L} p_{a,a'} E_{a'},
\end{equation}
where only finitely many of the coefficients $p_{a,a'}$ are non-zero, and $p_{a,a'}\in v \ZZ[v]$.
Expanding both sides of the equality \eqref{eq:b-bar-fixed} in the
basis $\{E_{a'}\}$, we rewrite it as the system of equations
\begin{equation}
\label{eq:p-condition}
p_{a,a'} - \overline{p_{a,a'}} =
r_{a,a'} + \sum_{a''} \overline{p_{a,a''}} r_{a'',a'}.
\end{equation}

Let us first show that the equations \eqref{eq:p-condition}
imply that $p_{a,a'}=0$ unless $a' \prec a$, hence the desired
element~$C_a$ must satisfy \eqref{eq:b-double-triangularity}.
Assume for the sake of contradiction that $p_{a,a'} \neq 0$ for
some $a' \not \prec a$, and choose $a'$ as some maximal element with this
property (this is possible since the sum in \eqref{eq:matrix-b-e}
is finite).
For this choice of $a'$, the right hand side of \eqref{eq:p-condition}
becomes~$0$. But since $p_{a,a'} \in v \ZZ[v]$, the condition
$p_{a,a'} - \overline{p_{a,a'}} = 0$ implies that $p_{a,a'} = 0$,
the desired contradiction.

To show the uniqueness of $p_{a,a'}$ for $a' \prec
a$, we proceed by induction on the maximal length of a chain
between $a'$ and $a$.
Thus, we can assume that the statement is already known for all
$p_{a,a''}$ appearing in the right hand side of \eqref{eq:p-condition}
since we must have $a' \prec a'' \prec a$ for the corresponding
term to be non-zero.
Now the uniqueness of $p_{a,a'}$ is obvious for the same reason as
above: any $p \in v \ZZ[v]$ is uniquely determined by $p - \overline p$.

To show the existence, we use the following obvious property: a
Laurent polynomial $f \in \ZZ[v, v^{-1}]$ can be written (uniquely) as $p -
\overline p$ for $p \in v \ZZ[v]$ if and only if $f + \overline f = 0$.
Indeed, we have $p = [f]_+$, where $[f]_+$ stands for the part of
the Laurent expansion of~$f$ that contains the positive powers of~$v$.
Thus, to check that \eqref{eq:p-condition} has a (unique) solution
for $p_{a,a'}$, it suffices to show that the the right hand side $f$ of \eqref{eq:p-condition}
satisfies $f + \overline f = 0$.
This is done by the following calculation using \eqref{eq:r-condition-twisted}
and our inductive assumption:
\begin{align*}
f + \overline f & = r_{a,a'} + \overline{r_{a,a'}} + \sum_{a''}
(\overline{p_{a,a''}} r_{a'',a'} +
p_{a,a''}\overline{r_{a'',a'}})\\
& = \sum_{a''} (- r_{a,a''}\overline{r_{a'',a'}} + \overline{p_{a,a''}} r_{a'',a'}
+ (\overline{p_{a,a''}} + r_{a,a''} + \sum_{a'''} \overline{p_{a,a'''}}
r_{a''',a''})\overline{r_{a'',a'}})\\
& = \sum_{a''} \overline{p_{a,a''}} (r_{a'',a'} + \overline{r_{a'',a'}}
+ \sum_{a'''} r_{a'',a'''}\overline{r_{a''',a'}}) = 0,
\end{align*}
as desired.

To complete the proof, it remains to show that our finiteness
assumption on the indexing poset~$L$ implies that, for a given $a \in L$,
only finitely many of the coefficients $p_{a,a'}$ are nonzero.
For $a' \preceq a$, let $c(a')$ denote the maximal length~$m$ of a
chain $a' = a_0 \prec a_1 \prec \cdots \prec a_m = a$ in~$L$.
Let $P_m(a)$ denote the set of elements~$a'$ such that $c(a') = m$
and $p_{a,a'} \neq 0$.
Since by our assumption, the function $c(a')$ is bounded from
above, it is enough to show that each $P_m(a)$ is finite.
Let $R(a)$ denote the (finite!) set consisting of the elements
$a' \in L$ such that $r_{a,a'} \neq 0$.
In view of \eqref{eq:p-condition}, we have
$$P_m(a) \subseteq R(a) \cup \bigcup_{k=1}^{m-1} \bigcup_{a'' \in P_k(a)}
R(a''),$$
hence the desired finiteness of $P_m(a)$ follows by induction on~$m$.
\end{proof}

\begin{remark}
Using the notation $[f]_+$ introduced in the course of the proof
of Theorem~\ref{th:Lusztig-lemma}, we can express the
coefficients $p(a,a')$ recursively by
\begin{equation}
\label{p-explicit}
p_{a,a'} =
[r_{a,a'} + \sum_{a''} \overline{p_{a,a''}} r_{a'',a'}]_+.
\end{equation}
\end{remark}

\end{document}